\documentclass[draft]{amsart}
\usepackage{a4}
\usepackage{latexsym}                  
\usepackage{exscale}                   
\usepackage[arrow,curve,matrix,cmtip,line]{xy}

\numberwithin{equation}{section}

\theoremstyle{plain}
\newtheorem{thm}{Theorem}[section]
\newtheorem{lem}[thm]{Lemma}
\newtheorem{prop}[thm]{Proposition}
\newtheorem{cor}[thm]{Corollary}

\NoComputerModernTips
\CompileMatrices
\xymatrixrowsep{1.5em}
\xymatrixcolsep{0.25em}
\newcommand{\krue}{1.1}                
\newcommand{\strich}{{\scriptscriptstyle /}}

\newcommand{\ane}[1]{`#1'}

\newcommand{\tst}{\textstyle}

\newcommand{\ssst}{\scriptscriptstyle}

\newcommand{\field}[1]{\mathbb{#1}}

\newcommand{\N}{\field{N}}
\newcommand{\No}{\N_0}
\newcommand{\R}{\field{R}}
\newcommand{\Rd}{\field{R}^d}
\newcommand{\Z}{\field{Z}}
\newcommand{\Zd}{\field{Z}^d}

\newcommand{\deff}{\stackrel{\rm def}{=}}

\DeclareMathOperator{\Id}{Id}
\DeclareMathOperator{\inter}{int}

\newcommand{\const}{K}
\newcommand{\conste}{L_1}
\newcommand{\constz}{L_2}
\newcommand{\grenze}{L}
\newcommand{\step}{\ell}
\newcommand{\decc}{L}
\newcommand{\sist}{S}              
\newcommand{\stwt}{u}              
\newcommand{\stv}{\underline{s}}   
\newcommand{\stl}{\ell'}           

\newcommand{\abs}[1]{\lvert#1\rvert}
\newcommand{\norm}[1]{\lVert#1\rVert}
\newcommand{\babs}[1]{\bigl\lvert#1\bigr\rvert}

\newcommand{\Babs}[1]{\Bigl\lvert#1\Bigr\rvert}

\newcommand{\normu}[1]{\norm{#1}_\infty}
\newcommand{\norme}[1]{\norm{#1}_1}

\newcommand{\normd}[1]{\norm{#1}_{\ssst\mathcal{D}}}
\newcommand{\normw}[1]{\norm{#1}_{\ssst\mathcal{W}}}

\newcommand{\AG}{G}
\newcommand{\AH}{H}
\newcommand{\ag}{g}
\newcommand{\ah}{h}

\newcommand{\ul}[1]{\underline{#1}}

\newcommand{\ulb}{\ul{b}}
\newcommand{\ulag}{\ul{\ag}}

\newcommand{\fAG}{\boldsymbol{\AG}}
\newcommand{\fag}{\boldsymbol{\ag}}

\newcommand{\fah}{\boldsymbol{\ah}}
\newcommand{\fB}{\mathbf{B}}

\newcommand{\fA}{\mathbf{A}}
\newcommand{\fa}{\mathbf{a}}

\newcommand{\tiaa}{\widetilde{A}}

\newcommand{\tiAG}{\widetilde{\AG}}
\newcommand{\tiag}{\widetilde{\ag}}
\newcommand{\tiah}{\widetilde{\ah}}

\newcommand{\tifAG}{\widetilde{\fAG}}
\newcommand{\tifag}{\widetilde{\fag}}
\newcommand{\tifah}{\widetilde{\fah}}

\newcommand{\ba}{\Bar{a}}
\newcommand{\bc}{\Bar{c}}
\newcommand{\baa}{\Bar{A}}
\newcommand{\bcc}{\Bar{C}}
\newcommand{\bbb}{\Bar{B}}
\newcommand{\bd}{\Bar{\delta}}
\newcommand{\bm}{\Bar{\mu}}

\newcommand{\sumx}[1][x]{\sum_{#1\in\Zd}}
\newcommand{\sumn}[1][1]{\sum_{n=#1}^\infty}
\newcommand{\summ}[1][1]{\sum_{m=#1}^n}
\newcommand{\summu}[1][1]{\sum_{m=#1}^\infty}
\newcommand{\summl}[1][l]{\sum_{m=1}^{#1}}
\newcommand{\sumj}[1][n+1]{\sum_{j=#1}^\infty}
\newcommand{\sums}[3][=n]{\sum_{\substack%
{\omega:#2\leadsto#3\\ \abs{\omega}#1}}}
\newcommand{\sumsx}[1][=n]{\sums[#1]{0}{x}}

\begin{document}
\title[A CLT for convolution equations]{A central limit theorem for
convolution equations and weakly self-avoiding walks}

\author{Erwin Bolthausen}
\author{Christine Ritzmann}
\date{\today}
\address{{\emph{Address of authors:}} Angewandte
Mathematik\\Universit\"at Z\"urich\\ Winterthurer Str.\,190\\CH-8057
Z\"urich-Irchel\\ Tel.\,xx41/1/635~58~49 (E.B.) and xx41/1/635~58~54
(C.R.)}
\email{eb@amath.unizh.ch\\chritz@amath.unizh.ch}

\begin{abstract}
The main result of this paper is a
general central limit theorem for distributions defined by certain renewal
type equations.  We apply this to weakly self-avoiding random walks. We
give good error estimates and Gaussian tail estimates which have not been
obtained by other methods.

We use the \ane{lace expansion} and at the same time develop
a new perspective on this method:  We work with a
fixed point argument directly in $\Zd$ without
using Laplace or Fourier transformation.
\end{abstract}
\maketitle
\setcounter{tocdepth}{1}
\tableofcontents
\section{Introduction and Results}
\label{chap-def}

\subsection{Introduction}

The standard simple random walk on the hypercubic lattice $\mathbb{Z}^{d}$ is
given by the uniform distribution on the set of nearest-neighbour paths
starting in $0,$ and of length $n.$ The law of the self-avoiding random walk
is simply the uniform distribution on the set of walks having \textit{no
self-intersections}. An interpolation between the strictly self-avoiding walk
and the standard random walk is the so-called weakly self-avoiding walk (also
known as \ane{Domb-Joyce model}). Here, self-intersec\-tions are not
completely forbidden, but penalized by a factor $1-\lambda$ for every
self-intersection, where $\lambda\in\left(  0,1\right)  $ is a parameter.
Not much is known rigorously for these self-avoiding random walks in
dimensions $d=2,3$ and $4.$

A quantity of basic interest is the so called connectivity $C_{n}\left(
x\right)  ,$ $x\in\mathbb{Z}^{d},$ $n\in\mathbb{N},$ which is obtained by
summing the weights of all paths from $0$ to $x$ of length $n$: In the random
walk case, the paths all get weight $1,$ in the strictly self-avoiding case,
only the paths without self-intersections get weight $1$, the others $0,$ and
in the weakly self-avoiding case, the weight is given in terms of the
number of self-intersections (and the parameter $\lambda)$ as indicated above.
After normalization, this defines the distribution of the end-point of the
walk. The first results in the case $d\geq5$ were obtained by Brydges and
Spencer \cite{BS} in the middle of the eighties. They introduced a
perturbative expansion technique, based on the so called lace expansion. With
this, Brydges and Spencer proved a central limit theorem for
the weakly self-avoiding walk in dimensions greater than four and for
small parameter
$\lambda$. The lace expansion is a kind of renewal equation for $C_{n}$ in
terms of the connectivities of the standard random walk of the following
type:%
\begin{equation}
C_{n}=C_{n-1}\ast 2d\, D+\sum_{m=1}^{n}\Pi_{m}\ast C_{n-m},\;n\geq
1.\label{Recursion}%
\end{equation}
Here $D$ is the usual nearest-neighbour distribution $D\left(  x\right)
=1/(2d)$ if
$\left|  x\right|  =1,$ and $0$ otherwise. The $\Pi_{m}$ reflect the
presence of the self-interaction of the paths. If they are not present,
then $C_{n}$ is of course just obtained by convoluting $2d\,D.$
(\ref{Recursion}) can easily be derived by a kind of inclusion-exclusion.
For the convenience of the reader, we give the derivation in Appendix B.
One of the delicacies of the problem is that the $\Pi_{m}$ are complicated
expressions which cannot be written down explicitly, but the key point is
that they can be estimated in terms of the
$C_{k},\,k\leq m.$ Similar expansions have now found widespread
applications in probability theory and mathematical physics, like in
percolation theory
\cite{HSc}, branched polymers and quite recently also for superprocess
approximations of interactive particle systems.

The strictly self-avoiding
walk was first treated by Slade \cite{Slade} for large dimensions,
where the large dimension serves in some sense as a small coupling
parameter. Finally, in a tremendous effort, Hara and Slade \cite{HSa,HSb}
were able to prove the diffusive behavior for $d\geq5$ for the
strictly self-avoiding case$.$ However, their argument is still
essentially perturbative and relies on computer-assisted estimates of a
number of constants which have to be smaller than $1.$

The earlier approaches to the lace expansion always depended on taking
(complex) Laplace transforms in time $\sum_{n}z^{n}C_{n}$, and then inverting
the transform. The latter step is notoriously difficult and
leads to a number of problems. In particular, it seems hard to obtain good
pointwise estimates for the connectivities in this way. In 1997,
van der Hofstad, den Hollander and Slade \cite{NLE} presented an
inductive approach to the lace expansion avoiding Laplace inversion. In
particular, these authors prove a local central limit theorem for
the so-called \ane{elastic} self-avoiding walk, a model in which the
penalty for self-inter\-sections decreases (only polynomially) in
time. Several further improvements have been obtained recently. In
\cite{GLE}, van der Hofstad and Slade generalize and simplify the
inductive approach. In a recent paper by Hara, van der Hofstad and Slade
\cite{HHS}, a general method for deriving information on the generating
two-point functions in $\mathbb{Z}^{d}$, i.e. expressions of the form
$\sum_{n}z^{n}C_{n}\left(  x\right)  $ at the critical point $z_{c}$ 
is developed.

In this paper we present a novel way to treat such problems. We first
modify (\ref{Recursion}) slightly. To explain the main idea, one has to
remark that the norming constants
$c_{n}\overset{\mathrm{def}}{=}\sum_{x}C_{n}\left( x\right)  $ behave in
first order exponentially in $n:$ $\lim_{n\rightarrow
\infty}\left(  1/n\right)  \log c_{n}$ exists$.$ It turns out that, in
leading order, the $\Pi_{n}$ have the same (exponential) behavior as the
$C_{n}.$ It is therefore natural to define
$B_{n}\overset{\mathrm{def}}{=}\Pi_{n}/c_{n},$ and rewrite
(\ref{Recursion}) as%
\begin{equation}
C_{n}=C_{n-1}\ast 2d\,D+\sum_{m=1}^{n}c_{m}B_{m}\ast C_{n-m},\;n\geq
1.\label{Recursion2}%
\end{equation}
The main result of our paper addresses the following problem: Given a sequence
$\left(  B_{m}\right)  _{m\geq1}$ of (possibly signed) distributions on
$\mathbb{Z}^{d}$, we regard (\ref{Recursion2}) as the defining equations for
the sequence $\left(  C_{n}\right)  _{n\geq0},$ of course with $C_{0}%
\overset{\mathrm{def}}{=}\delta_{0}.$ We prove that if the sequence $\left(
B_{m}\right)  $ has suitable smallness and decay properties, then $C_{n}%
/c_{n}$ satisfies a central limit theorem. It is important that we require
only polynomial decay properties of the $\left(  B_{m}\right)  .$ (In
fact, one crucial difficulty for analyzing the (weakly) self-avoiding
walks is that there is no mass gap between the $C$-sequence and the
$\Pi$-sequence). Although, the main result is completely independent of
the problem of self-avoiding walks, the version of the central limit
theorem we prove is tailored much with this application in mind. We come
as close as possible (we think) to a local central limit theorem,
including good Gaussian tail estimates. The (weakly) self-avoiding walk
can in fact not satisfy a local CLT in the strict sense.\footnote{The
walk starts at zero and will not forget this fact ever nor will the
memory weaken. The resulting error term is of size
$n^{-d/2}$, the same as the value of the approximating normal density at
zero.}.

The main idea of our approach is to regard (\ref{Recursion2}) as a fixed-point
equation of an operator acting on sequences of distributions. There are of
course many ways to do that, and the problem is to find the \ane{right}
way. We prove a central limit theorem by the following reasoning: The
operator we define has two crucial properties: First, it has contraction
properties in suitable spaces, so that a fixed point exists by the Banach
fixed point theorem. This fixed point has to satisfy (\ref{Recursion2}),
and as this has a unique solution (under appropriate conditions), we know
that the fixed point of our operator is \textit{the }solution. Secondly,
we prove that the operator leaves the property of being
\ane{asymptotically Gaussian} untouched: If a sequence is asymptotically
Gaussian (in a sense which has to be quantified, and after
normalization), then the same is true for the sequence transformed by the
operator. If we start with a sequence being already asymptotically
Gaussian (for instance just taking the discretized Gaussian
distribution), then also the fixed point is asymptotically Gaussian, and
therefore, the solution of (\ref{Recursion2}) satisfies a central limit
theorem (after normalization). This approach reveals in a transparent way
what is behind the central limit theorem: The property of being
asymptotically Gaussian is stable enough not to be destroyed by the
presence of the $B$-sequence, provided this enjoys the appropriate
properties.

As remarked above, the result we prove is tailored much to make the
application to weakly self-avoiding walks very easy. The point is that the
$B^{\prime}$s can be estimated in terms of the $C^{\prime}$s$.$ These
estimates are most naturally in $x$-space and not in Fourier-space. For this
reason, we work all the time in $x$-space, except for deriving several crucial
properties of convolutions of $D.$ By these estimates, we prove by a simple
\ane{circular} argument, that the \ane{true} $B^{\prime}$s coming from
the weakly self-avoiding walk do in fact satisfy the conditions of our
central limit theorem.

We believe that the method can be extended in various ways to treat many
other problems as well, but here we restrict ourselves to the
application to weakly self-avoiding walks. Even for this well studied
problem, the results we obtain are sharper than those obtained by other
methods.

\subsection{The Objects of the Weakly Self-Avoiding Walk}

We begin by introducing the connectivities $C_n$ for the weakly
self-avoiding walk. For
$x\in\Zd$ we set
$C_0(x)\deff \delta_{0,x}$, and for $n\ge 1$ and $\beta\ge 0$, we
define
\begin{equation}                                           \label{Cndef}
C_n(x) \deff \sumsx e^{-\beta\sum_{0\le s<t\le n}U_{st}(\omega)}
= \sumsx \prod_{0\le s<t\le n} (1-\lambda U_{st}(\omega)),
\end{equation}
where the sum is over all $n$-step simple random walk paths $\omega$
from $0$ to
$x$, while
\begin{align}                                     \notag
U_{st}(\omega) &\deff \begin{cases} 1,&\text{if }\omega(s)=\omega(t),\\
                                             0,&\text{if not,}\end{cases}\\
\intertext{and}
\lambda &\deff 1-e^{-\beta}.\notag
\end{align}
Here, loops of the walk are penalized with a factor $1-\lambda$.
We have the fully self-avoiding walk for $\lambda=1$, while for
$\lambda=0$ it is the simple random walk. We will
usually suppress the
$\lambda$-dependence in our notation.

The lace expansion is a renewal type equation for the connectivities.
It involves a function $\Pi_m$, defined in \eqref{Pimdef}, which we will call
\ane{lace function}. The expansion is the recursion formula stating that for
all $x\in \Zd$ we have
\begin{equation}                                          \label{recursion}
C_n(x) = 2d(D*C_{n-1})(x) + \sum_{m=2}^n
\Pi_m*C_{n-m}(x),
\end{equation}
where $D$ denotes the law of one step of a simple random walk, that is,
\begin{equation*}
D(y)\deff \begin{cases}
                               1/(2d),&\text{if }\norme{y}=1,\\
                               0,&\text{else.}
                  \end{cases}
\end{equation*}
The star $*$ always refers  to the discrete folding of two measures on $\Zd$,
that is $\AG*\AH(x)\deff \sum_{y\in\Zd}\AG(y)\AH(x-y)$.
Denote by $c_n$ the total mass of $C_n$, that is,
\begin{align}                                    \label{cndef}
c_n &\deff \sumx C_n(x).
\end{align}
(We will always denote measures by capital letters and the corresponding total
mass by the appropriate lower case letter.)
Furthermore we will always write $\varphi_\eta$ for
the d-dimensional normal density with
covariance matrix $\eta\cdot \Id_d$, that is,
\begin{align}                                     \notag
\varphi_\eta(x)
\deff (2\pi\eta)^{-d/2} \;\exp \bigl(\tst-\frac{x^2}{2\eta}\bigr).
\end{align}

The general results proved in the sections \ref{chap-const} and
\ref{chap-local} can be applied to the weakly self-avoiding walk ($\lambda$
small enough) in dimensions above four. This is stated in the following
theorem. The first part recovers the result of Brydges and Spencer
\cite{BS}. The second part presents a pointwise comparison of the
probabilities $C_n(x)/c_n$ with a normal
density and gives strong error estimates.
\begin{thm}                                             \label{thm-main}
Let the dimension $d$ be greater or equal to five.
Then there exists a
$\lambda_0>0$ such that for all
$0\le\lambda\le \lambda_0$ and $n\in\N$,
\begin{align}                                           \label{main1}
c_n &=\alpha \mu^n \bigl(1+ O(n^{-1/2})\bigr).\\
\intertext{For $x\in\Zd$ and $n\in\N$ such that $n -\norme{x}$ is even,%
\footnotemark \, we have}
\Babs{\frac{C_n(x)}{c_n} - 2\varphi_{\delta n}(x)}
&\le \const  \Bigl[n^{-1/2} \varphi_{n \nu}(x) + n^{-d/2} \sum_{j=1}^{n/2} j
\varphi_{j
\nu}(x)\Bigr].                                              \label{main2}
\end{align}
The constants $\alpha$, $\mu$ and $\delta$ are positive and depend on
$\lambda$ and $d$, whereas $\nu$ and $\const $ only depend on the dimension.
\end{thm}

\footnotetext{Because of the two-periodicity of nearest neighbour
walks, $C_n(x)$
equals zero when $\norme{x}$ and $n$ don't have the same parity.
The factor $2$ for the normal density in
\eqref{main2} arises for the same reason -- we have to double the density
values to approximate a discrete probability measure on the two-periodic
sublattices of $\Zd$.}

Note that for \ane{large} $x$ the right-hand side of \eqref{main2} is
bounded by
a multiple of
$n^{-1/2}\varphi_{n\nu}(x)$ alone. For $\abs{x}\le O(n^{1/2})$ the leading term
in the error bound is the second part, which is of order $O(n^{-d/2})$ near
zero. The constants
$\alpha$,
$\mu$ and
$\delta$ are identified in terms of $c_n$ and $\pi_m$ at the end of section
       \ref{chap-saw}.

\subsection{Strategy and Motivation of the Approach}

We prove the theorem by splitting the problem into two parts:
In the first step we show the existence
of the connective constant
$\mu$ and the diffusion constant
$\delta$ and we determine their exact form. Given these parameters, we can
write down the proper normal density to approximate the distributions
themselves. So in the second step we only have to estimate the error of the
approximation.

Several difficulties arise in this approach.
The constants we want to determine in advance are
given only implicitely as series in terms of $c_n$ and $\pi_n$ (the total
mass constants of the \ane{lace
functions} $\Pi_n$). In particular,  we need the total mass values of the
$n$-step weakly self-avoiding walk {\em for each natural $n$}.
So at
first it seems impossible to determine the constants without knowing
quite a lot about the distributions themselves.

We deal with this difficulty by treating the whole sequence
of mass constants as one object instead of studying the
constants separately for each $n$.
The idea of working in this kind of sequence spaces was taken from Gr\"ubel
\cite{G}.
We can define an operator on some appropriate sequence space such that
the sequence satisfying the recursion formula of the lace expansion is a
fixed point of this operator. Once we have chosen space and operator properly,
the Banach fixed point theorem yields the desired properties of the sequence.

The second step is organized as follows. To start we define a discrete random
vector, whose variance is given by the diffusion constant from the first
part. Then we consider the measures of the weakly self-avoiding walk as
perturbations of the simple random walk whose single steps are given by the
random vector defined above. We again use a fixed point argument to
control the errors of the approximation, this time working in some
sequence space of measures rather than of real numbers.

To give the fixed point arguments, the sequences we investigate shouldn't grow
exponentially. So first we cancel the exponential growth of the mass
constants for both the connectivities and the lace functions:
For $m\ge 2$d we define the function $B_m$ on $\Zd$ by
\begin{equation}                                        \label{Bmdef}
B_m(x) \deff \frac{\Pi_m(x)}{\lambda c_m}.
\end{equation}
Let $b_m$ denote the total mass of $B_m$.
Inserting \eqref{Bmdef} into \eqref{recursion}, we obtain the following
recursive identities (from now on suppressing $x$ in the notation):
\begin{align*}
C_n =& 2d D*C_{n-1} +\lambda \sum_{m=2}^n c_m B_m*C_{n-m}\quad\text{and}\\
c_n =& 2d c_{n-1} +\lambda \sum_{m=2}^n c_m b_m c_{n-m}.
\end{align*}

Now suppose that $C_n$ grows exponentially, that is, $C_n$ equals $\mu^n A_n$
with some
$\mu>0$ such that $a_n=\sumx A_n(x)$ tends to some $\alpha\neq0$ when
$n$ tends to
infinity. This $\mu$ is called {\em connective constant}.
The identities above lead to the following equations:
\begin{align}
A_n =& 2d\mu^{-1}D*A_{n-1} +\lambda \sum_{m=2}^n a_m B_m*A_{n-m},
\label{def-An}\\
a_n =& 2d\mu^{-1}a_{n-1} +\lambda \sum_{m=2}^n a_m b_m a_{n-m}.
\label{def-an}
\end{align}
If the sequence $(a_n)$ is converging to a limit $\alpha>0$ and
if the $b_n$ decay fast enough, we can let $n$
tend to infinity in \eqref{def-an} to find
\begin{equation*}
\alpha = 2d\mu^{-1}\alpha +\lambda \sum_{m=2}^\infty a_m b_m \alpha,
\end{equation*}
and therefore we have for the {\em connective constant} $\mu$,
\begin{equation}                                          \label{mu}
2d\mu^{-1} = 1 - \lambda \sum_{m=2}^\infty a_m b_m.
\end{equation}
By substituting this into \eqref{def-an} we obtain
\begin{equation}                                   \notag 
a_n = a_{n-1} -\lambda \left(\sum_{m=2}^\infty a_m b_m a_{n-1} - \sum_{m=2}^n
a_m b_m a_{n-m}\right).
\end{equation}
This derivation in mind, we will prove the existence and uniqueness of such
a sequence by a Banach fixed point argument in section \ref{chap-const},
under the hypothesis that the $b_m$ decay fast enough.

Given the sequence $(a_n)$ and specific
pointwise estimates of $B_m(x)$, we will obtain pointwise approximations
for $A_n(x)$ in dimension $d\ge 5$ in section \ref{chap-local}. This central
limit theorem with error estimates is our main result. The
theorem is proven for sequences satisfying a slightly
generalized version of \eqref{def-An}.

In section \ref{chap-saw} we prove the right behavior of the lace
expansion terms, insuring that we can indeed apply the fixed point arguments to
the weakly self-avoiding walk.
\section{Determining the Mass Constants}\label{chap-const}

\subsection{Existence and Uniqueness}

Let $(b_m)_{m\in\N}$ be a realvalued sequence with
\begin{align}                                       \notag 
\beta\deff \summu m\,\abs{b_m} < \infty.
\end{align}

In this section we will prove the following result:
\begin{prop}                                               \label{prop-const}
There is a $\lambda_0 = \lambda_0(\beta)>0$ such that for all
$\lambda\le\lambda_0$ there exists a unique sequence
$(a_n)_{n\in\No}$ with
\begin{align}                                                \label{aconst}
a_0=1 \quad\text{and}\quad
a_n =\bigl(1 -\lambda \summu a_m b_m \bigr) a_{n-1}
+ \lambda\summ a_m b_m a_{n-m},
\end{align}
such that $\sum_{n=1}^\infty \abs{a_n-a_{n-1}}<\infty$.
\end{prop}

We will prove this proposition with a fixed point argument, but first we
introduce some notation. Let $(l_\infty,\normu{.})$ be the Banach space of
bounded real valued  sequences $\fag=(\ag_n)_{n\in\No}$ with the supremum norm.
The difference operator $\Delta:\R^{\No} \to \R^{\No}$ is given by
\begin{equation}                                  \notag 
(\Delta \ag)_0 \deff  \ag_0        \qquad\text{and}\qquad
(\Delta \ag)_n \deff  \ag_n - \ag_{n-1} \qquad\text{for }n\in\N.
\end{equation}

For $\fag = (\ag_n)_{n\in\No}$ with $\sumn[0]
\abs{(\Delta \ag)_n}<\infty$ define the norm
\begin{align}                                     \notag 
\normd{\fag}&\deff \sumn[0] \abs{(\Delta \ag)_n}.
\end{align}
Furthermore, define the operator $\sim$ on sequences by
\begin{align}                                                \label{simdef}
{\tiag}_0&\deff \ag_0  \qquad\text{and}\notag\\[1mm]
\tiag_n
&\deff \tiag_{n-1} - \lambda
\Bigl[\summ \ag_m b_m(\ag_{n-1}-\ag_{n-m}) + \ag_{n-1} \sumj
\ag_j b_j \Bigr].
\end{align}

We will apply the Banach fixed point theorem to the operator $\sim$. In the
three following lemmas  we prove that
the necessary conditions are fulfilled.

\begin{lem}                                             \label{lem:op-eins}
Let $\fag\in l_\infty$ with $\normd{\fag}<\infty$. Then we also have
$\normd{\tifag}<\infty$.\end{lem}
\begin{proof}
Let $\normd{\fag}<\infty$.
We have to show that $\sumn[0]\abs{(\Delta \tiag)_n}$
is finite. First notice that
\begin{equation}                                          \label{normu}
\normu{\fag}= \sup_{n\in\No}\, \abs{\ag_n}
=\sup_{n\in\No} \,\abs{\sum_{k=0}^n (\Delta \ag)_k}
\le \normd{\ag}.
\end{equation}
          From $\eqref{simdef}$ we have
\begin{align}                                             \label{Delta_Ft_n}
\sumn[0]\abs{(\Delta \tiag)_n}
&= \abs{\ag_0}+\sumn\abs{\tiag_n-\tiag_{n-1}}
\notag\\
&= \abs{\ag_0}+\lambda\sumn\Babs{\summ \ag_m b_m
(\underbrace{\ag_{n-1}-\ag_{n-m}}_{=\sum_{l=1}^{m-1} (\Delta \ag)_{n-l}}) +
\sumj \ag_j b_j \ag_{n-1}}\notag\\
&\le \abs{\ag_0}+\lambda\Bigl[ \normu{\fag}\sumn\sum_{l=1}^{n-1}\,\,
\abs{(\Delta
\ag)_{n-l}}\!\!\!\underbrace{\summ[l+1]\abs{b_m}}_{\le\summu[l+1]\abs{b_m}}
\!\!+ \normu{\fag}^2
\underbrace{\sumn\sumj\abs{b_j}}_{\le \sum_{j=2}^\infty j\abs{b_j}\le \beta}
\Bigr]\notag\\
&\le \abs{\ag_0}+\lambda \Bigl[\normu{\fag} \sum_{l=1}^{\infty}
\summu[l+1]\abs{b_m}\sumn\abs{\Delta\ag_n}
+ \normu{\fag}^2 \beta\Bigr]\notag\\
&\le \abs{\ag_0}+2\lambda\beta \normd{\fag}^2 <\infty,
\end{align}
where we used \eqref{normu} in the last line.
\end{proof}

\begin{lem}                                            \label{lem:contract}
Let $\mathcal{D}_\grenze \deff\{\ag\in l_\infty: \ag_0=1 \text{ and
}\normd{\fag}\le \grenze\}$, where
$\grenze$ is a constant greater than or equal to $3/2$. Then for all
$\lambda\le 1/(6\beta \grenze)$ the operator $\sim$ is a contraction
with respect
to $\normd{.}$ on $\mathcal{D}_\grenze $.
\end{lem}

Note that the value $3/2$ in the lemma is chosen to keep
the constants simple. An analogous statement holds as long as $\grenze$ is
bounded away from one.

\begin{proof}
We have to show:
\begin{itemize}
\item[(i)]
$\fag\in \mathcal{D}_\grenze \;\Rightarrow \;\tifag\in
\mathcal{D}_\grenze $
\item[(ii)]
There exists some $\kappa<1$
such that$\normd{\tifag-\tifah}\le
\kappa\normd{\fag-\fah}$ for all $\fag,\fah\in \mathcal{D}_\grenze $.
\end{itemize}

To see (i),
let $\fag\in\mathcal{D}_\grenze$ be given. We know $\tiag_0=\ag_0=1$.
According to \eqref{Delta_Ft_n} we have
\begin{equation*}
\normd{\tifag}
\le 1+2 \lambda\beta  \grenze^2 \le \frac{2}{3} \grenze + \frac{1}{3}\grenze,
\end{equation*}
whenever $\lambda\le\frac{1}{6\beta \grenze}$.

To see (ii),
take $\fag,\fah\in \mathcal{D}_\grenze $. Since $\ag_0$ equals $\ah_0$,
we have
$\normd{\tifag-\tifah}=\sumn \abs{(\Delta (\tiag-\tiah))_n}$.
We have (from \eqref{simdef}):
\begin{align}                                          
\sumn\abs{(\Delta \tiag)_n-(\Delta \tiah)_n}
&= \lambda\,\sumn\,
\biggl\lvert\,\summ (\ag_m - \ah_m)b_m (\ag_{n-1}-\ag_{n-m}) \notag\\
&\qquad\qquad\;+ \summ
\ah_m b_m \bigl[\ag_{n-1}-\ah_{n-1}-(\ag_{n-m}-\ah_{n-m})\bigr]\notag\\
&\qquad\qquad\;+ \sumj b_j
\bigl[\ag_j (\ag_{n-1}-\ah_{n-1})
                             + (\ag_j -\ah_j)\ah_{n-1}\bigr]
\biggl\rvert\notag
\end{align}
We can estimate the absolute values of the three summands individually. The
first one can be treated analogously to \eqref{Delta_Ft_n},
\begin{align}                               \notag 
\lambda\sumn\summ \abs{b_m}\, \abs{\ag_m -\ah_m}\,
\abs{\ag_{n-1}-\ag_{n-m}}
&\le \lambda\beta\normu{\fag-\fah} \normd{\fag}.
\end{align}
Very similarly we obtain for the second one
\begin{align}                                  \notag 
\lambda\sumn\summ \abs{b_m}\,
\abs{\ah_m}\,
\abs{\ag_{n-1}-\ah_{n-1}-(\ag_{n-m}-\ah_{n-m})}
&\le \lambda\beta \normu{\fah}\normd{\fag-\fah},
\end{align}
and for the third
\begin{align}                                    \notag 
           \lambda\sumn\sumj \abs{b_j}
\,\abs{\ag_j (\ag_{n-1}-\ah_{n-1})
                             + (\ag_j -\ah_j)\ah_{n-1}}
&\le \lambda\beta \normu{\fag-\fah}(\normu{\fag}+\normu{\fah}).
\end{align}
Since both $\fag$ and $\fah$ are in $\mathcal{D}_\grenze $ and
$\lambda\le\frac{1}{6\beta \grenze}$,
this yields
\begin{flalign}                                  \notag 
&&\normd{\tifag-\tifah}
&\le 4\lambda\beta \grenze \normd{\fag-\fah}\le \frac{2}{3}\normd{\fag-\fah}.
&&\qed\end{flalign}
\renewcommand{\qedsymbol}{}
\end{proof}

\begin{proof}[Proof of Proposition \ref{prop-const}]
The elements of $l_\infty$ with finite $\normd{.}$-norm form a Banach
space with
this norm, and $\mathcal{D}_\grenze$ is a closed subset of this  space.

Thus, using Lemma \ref{lem:op-eins} and Lemma \ref{lem:contract}, the
Banach fixed point theorem yields for small enough $\lambda$ the existence
and uniqueness of an element $\fa\in \mathcal{D}_\grenze $ with
$\widetilde{\fa} = \fa$. Furthermore, the repeated iteration of $\sim$
with starting point
$(1,1,1,\dots)$ converges to $\fa$. As long as $\grenze\ge 3/2$, the value
of
$\grenze$ has an influence only on the upper bounds for $\lambda$.
This proves the proposition.
\end{proof}

\subsection{Limit and Convergence Speed}

Now we investigate the limit and the convergence behavior of this
\ane{fixed sequence}
in a more particular setting. By choosing $\grenze= 3/2$ we obtain for all
$\lambda\le 1/(9\beta)$ a sequence $\fa$ with $a_0=1$,
$\sumn \abs{(\Delta a)_n}\le 1/2$ and for all $n\in\N$
\begin{align*}
a_n &= \stwt\mu^{-1} a_{n-1} + \lambda \summ a_m b_m a_{n-m},
\end{align*}
where $\stwt\mu^{-1} = 1 - \lambda \summu a_m b_m$ as in
\eqref{mu} for $\stwt=2d$. Since $\fa$ is bounded and
$\sumn[2]\abs{b_n}<\infty$, $\stwt\mu^{-1}$ is
finite. Note also that for all $n\in\No$ we have
\begin{equation}                                \notag 
1/2\le a_n \le 3/2.
\end{equation}

We now want to investigate the limiting value $\alpha
=\lim_{n\to\infty} a_n$. Since the difference sequence of $\fa$ is
absolutely summable, $\alpha$ exists, and we have
\begin{align*}
\alpha &= \lim_{n\to\infty} a_n = \lim_{n\to\infty} \summ[0](\Delta a)_m.
\end{align*}
Now consider for fixed $n\in\N$ (recall \eqref{aconst}):
\begin{align}
a_n &=  1 + \sum_{k=1}^n (\Delta a)_k\notag\\
&= 1- \lambda \sum_{k=1}^n \Bigl[\summu a_m b_m a_{k-1}
- \sum_{m=1}^k a_m b_m a_{k-m}\Bigr]\notag\\
&= 1- \lambda \summu a_m b_m \sum_{k=1}^n a_{k-1}
+ \lambda\summ a_m b_m \sum_{k=m}^n a_{k-m}\notag\\
&= 1- \lambda
\underbrace{\sum_{m=n+1}^\infty a_m b_m \sum_{k=1}^n a_{k-1}}_{=:F_1}
- \lambda\summ a_m b_m \sum_{l=1}^{m-1}
           \underbrace{a_{n-l}}_{=\alpha-(\alpha-a_{n-l})}\notag\\
&= 1- \lambda F_1
- \lambda\summ a_m b_m (m-1)\alpha
+ \lambda
\underbrace{\summ a_m b_m\sum_{l=1}^{m-1}\sum_{k=n-l+1}^\infty
(\Delta a)_{k}}
_{=:F_2},                                \label{alpha}
\end{align}
where
\begin{align*}
\abs{F_1} &= \Babs{\sum_{m=n+1}^\infty a_m b_m \sum_{k=1}^n a_{k-1}}
\le \sum_{m=n+1}^\infty \grenze \abs{b_m} n \grenze
\le \grenze^2 \sum_{m=n+1}^\infty m \abs{b_m}
\stackrel{\ssst n\to\infty}{\longrightarrow} 0\notag
\end{align*}
and
\begin{align*}
\abs{F_2} &= \Babs{\summ a_m b_m\sum_{l=1}^{m-1}\sum_{k=n-l+1}^\infty
(\Delta a)_{k}}
\le \sum_{l=1}^{n-1}\sum_{m=l+1}^\infty \grenze \abs{b_m}\sum_{k=n-l+1}^\infty
\abs{(\Delta a)_{k}}\notag\\
&\le \grenze \underbrace{\sum_{l=1}^{n/2}\sum_{m=l+1}^\infty
\abs{b_m}}_{\le \beta}
\sum_{k=n/2}^\infty \abs{(\Delta a)_{k}}
+ \grenze \sum_{l=n/2}^{n-1}\sum_{m=l+1}^\infty \abs{b_m}
\underbrace{\sum_{k=1}^\infty \abs{(\Delta a)_{k}}}_{\le \normd{\fa}}
\stackrel{\ssst n\to\infty}{\longrightarrow} 0.\notag
\end{align*}

Letting $n$ tend to infinity in \eqref{alpha}, we obtain
\begin{equation*}
\alpha = 1- \lambda \summu (m-1) a_m b_m \alpha,
\end{equation*}
which yields
\begin{equation}                                         \label{alpha1}
\alpha^{-1} = 1+ \lambda \summu (m-1)\, a_m b_m
= \stwt\mu^{-1} + \lambda \summu m\, a_m b_m.
\end{equation}

In case we know the rate of decay of the $b_m$, we can determine the
speed of the
convergence $a_n\longrightarrow \alpha$ more precisely. The following
corollary states a result that we will need in the next section.
\begin{cor}                                             \label{cor:speed}
If there exist positive constants $\varepsilon$ and $\beta'$ such that
\begin{equation*}
\abs{b_m}\le \beta' m^{-2-\varepsilon}\quad\text{for all }m\in\N,
\end{equation*}
then we get a decay of order $n^{-1-\varepsilon}$ for the difference sequence
$\Delta \fa$. More precisely we have
\begin{equation*}
\abs{(\Delta a)_n}\le \lambda \beta' \const  n^{-1-\varepsilon}\quad
\text{for all }n\in\N,
\end{equation*}
where $\const $ is a positive constant not depending on $\lambda$ or $\beta'$.
In particular we have another constant $\const $ such that
\begin{equation}                             \notag 
\abs{\alpha - a_n}\le \lambda \beta' \const  n^{-\varepsilon}\quad
\text{for all }n\in\N.
\end{equation}
\end{cor}
\begin{proof}
Using \eqref{alpha}, both estimates can be easily seen by induction.
\end{proof}
\section{Local Estimates in High Dimensions}
\label{chap-local}

We now turn to estimates not only for the normalization constants,
but for the measures on $\Zd$ themselves. In particular, we are interested in
the following question: If we consider the measures as
perturbations of the distribution of a sum of independent,
identically distributed
random vectors, then how big is the pointwise difference between the measure
and the appropriate normal density? Having in mind the high-dimensional
self-avoiding walk, we will not  expect a proper local central limit theorem to
hold: There will always be correction terms of order
$n^{-d/2}$ near zero, since zero is the starting point of the walk.
What we obtain is Gaussian decay for the perturbative errors on the
whole $\Zd$, improved by a factor of $n^{-1/2}$ for large $x$.

In this section, we will show the pointwise estimates in a more general
context. Supposing some specific pointwise bounds for the
distributions $B_n$, it is possible to show local estimates for the measures
$A_n$ in five or more dimensions. In the next section we will
show that the lace functions in the weakly self-avoiding walk case have
the necessary properties.

We begin by introducing some notations:
Let the space $\mathcal{M}$ be defined as the set of the symmetric,
rotationally invariant, signed real valued measures on $\Zd$
with existing \ane{variance}, that is,
\begin{multline}                                      \notag 
\mathcal{M}:=\{\AG:\Zd\to\R \text{ such that }\sumx \abs{\AG(x)}<\infty
\text{ and }\sumx x^2 \abs{\AG(x)}<\infty;\\
\AG\text{ symmetric in each coordinate and rotationally invariant}\}.
\end{multline}
Here and hereafter, we denote by $x^2$ the square of
the euclidean norm of $x\in \Rd$.
For $\AG\in\mathcal{M}$ define
\begin{equation}
\ag\deff\sumx \AG(x) \qquad\text{and}\qquad
\ulag\deff \sumx x^2\,\AG(x).                \notag 
\end{equation}
By $\mathcal{S}$ we denote the  space of sequences with elements in
$\mathcal{M}$, that is,
\begin{align}                                   \notag 
\mathcal{S}:=  \bigl\{\fAG= (\AG_n)_{n\in\No}: \AG_n\in\mathcal{M}\,
\text{ for all }n\in\No\bigr\}.
\end{align}

From the first section recall the equation \eqref{def-An},
\begin{align*}
A_n =& 2d\mu^{-1}D*A_{n-1} +\lambda \sum_{m=2}^n a_m B_m*A_{n-m}.
\end{align*}

Here we will investigate a slightly
different sequence $(A_n)$ with
\begin{align}                                           \label{A-unwohl}
\qquad A_n=
\stwt\mu^{-1}\sist *A_{n-1} +\lambda \summ a_m B_m*A_{n-m},
\end{align}
where $\stwt>0$ is a fixed constant ($\stwt=2d$ in
\eqref{def-An}). We suppose $S\in{\mathcal M}$ is a non-degenerate%
\footnote{that is, the probability $S(0)<1$.}
and aperiodic%
\footnote{Consider a probability measure $\sist$ and the random walk
defined as sum of independent steps, where each step is distributed
according to $\sist$. Now take the greatest common divisor of all times
$n$, at which the probability of staying at zero is not zero. This
divisor is called {\em period} of the walk. If the period equals one,
the random walk, and hence $\sist$ itself, is called {\em aperiodic}.}
probability measure of bounded range. More precisely, we have a
constant $\stl\ge1$ such that $\sist(x)=0$ for all $x$ with
$\abs{x}>\stl$. For simplicity we also assume that the
\ane{variance} $\stv$ is greater or equal to one. This condition is
not necessary,
but convenient, since many constants depend on the lower bound of the involved
variance.

Note that the distribution $D$ from \eqref{def-An} is {\em not}
aperiodic: For technical reasons we will give the proof for aperiodic
measures only, and afterwards discuss the two-periodic case to which the
self avoiding walk belongs.

For the whole section, the dimension $d$ is greater or equal to five.
$\const $ denotes a positive constant depending on $d$ and $\stl$ only.
The value of
$\const $ may change from line to line, whereas $\nu\ge 1/(2d)$ is an
adjustable
parameter and will be determined later.

Now let $\fB\in\mathcal{S}$ be a sequence with
$B_0\equiv 0$ and with
\begin{align}                                              \label{Bmlokprop}
\abs{B_m(x)} &\le \const  m^{-d/2}\sum_{k=1}^{m/2}k^{1-d/2}\varphi_{k \nu}(x)
\end{align}
uniformly for all $m\ge 1$ and $x\in\Zd$. We abbreviate the
notation by writing
\begin{equation*}
\psi_m (x) \deff m^{-d/2}\sum_{k=1}^{m/2}k^{1-d/2}\varphi_{k
\nu}(x)
\end{equation*}
and define
\begin{align*}                                           
\beta_\nu &= \sup_{m,x} \frac{\abs{B_m(x)}}{\psi_m(x)}.
\end{align*}

We have (recall $d\ge 5$)
\begin{align}                                           \label{bmdecay}
\abs{b_m} &\le \sumx \abs{B(x)}
\le \beta_\nu m^{-d/2}\sum_{k=1}^{m/2}k^{-3/2}
\sumx\varphi_{k \nu}(x) \le \const  \beta_\nu m^{-d/2}\quad\text{and}\\
\abs{\ulb_m}&\le\sumx x^2\,\abs{B(x)}
\le \beta_\nu
m^{-d/2}\sum_{k=1}^{m/2}k^{-3/2}
\sumx x^2\varphi_{k \nu}(x)\le \const \nu \beta_\nu m^{-(d-1)/2}.\notag
\end{align}
Here we used that as long as $\eta\ge 1/2d$, we have
\begin{equation*}
\sumx \varphi_{\eta}(x)\le \const \qquad\text{and}\qquad
\sumx x^2\varphi_{\eta}(x) \le \const  \eta,
\end{equation*}
which is stated and proved in Lemma \ref{lem-discsum} in the appendix. In
particular, since
$d\ge 5$, the condition on
$\beta$ from the last chapter is fulfilled, that is, $\beta=\summu
m\abs{b_m}\le \const \beta_\nu$.
In fact, the $b_n$ fulfill even the stronger condition from Corollary
\ref{cor:speed} for $\varepsilon=1/2$ and $\beta'=\sup_{m}m^{d/2}\abs{b_m}$.

For $\lambda$ small enough, Proposition \ref{prop-const} and Corollary
\ref{cor:speed} now yield the existence of a unique sequence
$(a_n)_{n\in\No}$ with
$a_0=1$ and
\begin{align*}
a_n &= \stwt\mu^{-1} a_{n-1} + \lambda \summ a_m b_m a_{n-m},
\end{align*}
where $\stwt \mu^{-1} = 1 - \lambda \summu b_m a_m$ and $(\Delta a)_n =
O(n^{-3/2})$.
Using these $a_n$, we are now in the situation to define the
sequence
$(A_n)$ of signed measures on $\Zd$ properly by (see \eqref{A-unwohl})
\begin{align}                                                 \label{defA}
A_0&\deff\delta_{0} \qquad\text{and} \qquad
A_n\deff \stwt\mu^{-1}\sist *A_{n-1} +\lambda \summ a_m B_m*A_{n-m}.
\end{align}

The reader might worry that there is a problem caused by the
ambiguous use of $a_n$, which denotes the $n$th term of the given fixed point
sequence on one hand and the total mass of
$A_n$ on the other hand. But by summing  up \eqref{defA} over $x\in\Zd$, we
see immediately that the normalization constant of $A_n$ in fact {\em is}
the given $a_n$.

We will use the following abbreviations:
\begin{align*}
\rho &\deff\summu a_m \,b_m,&
\sigma &\deff\summu (m-1) a_m \,b_m,&&\text{and}&
\tau &\deff \summu a_m \,\ulb_m.
\end{align*}

With these notations we have
\begin{align*}
\stwt\mu^{-1} &=1-\lambda\rho\qquad\text{and}\qquad
\alpha^{-1} =1+\lambda\sigma,
\end{align*}
where $\alpha$ is the limit of the sequence $(a_n)$ (see \eqref{alpha1}).
In addition we know that
$1/2\le a_n \le 3/2$ for all $n\in\No$, and that
$\alpha-a_n=O(n^{-1/2})$, which results from Corollary \ref{cor:speed}.

The key parameter of the following approximation is the constant
\begin{align}                                              \label{deltadef}
\delta \deff \frac{\stv(1-\lambda\rho) +\lambda\tau}{d(1+\lambda\sigma)},
\end{align}
which will turn out to be the right diffusion constant for the asymptotic
probability law $A_n(x/\sqrt{n})/a_n $. We always assume $\lambda$ to be small
enough to ensure that
\begin{equation}                                            \label{deltaest}
\stv/2\le d \delta \le 2\stv.
\end{equation}

Since $\stv\ge1$, we have in particular $\delta\ge 1/(2d)$. Now we
can state the
main result of this paper:
\begin{thm}[Local Estimates, aperiodic case]              \label{thm-local}
The sequence $(A_n)_{n\in\No}$ defined in \eqref{defA}, has
the following property: There exist $\lambda_0>0$ such that for all
$\lambda\in(0,\lambda_0)$ and for all
$x\in\Zd$:
\begin{align}                                             \label{local}
\abs{A_n(x) - a_n\varphi_{n\delta}(x)} &\le \const \Bigl[
n^{-1/2}\varphi_{n \nu}(x) + n^{-d/2}\sum_{j=1}^{n/2} j\,\varphi_{j
\nu}(x)\Bigr],
\end{align}
where $\const$ and $\nu$ are positive constants depending on $d$ and
$\stl$ only. (In particular they do not depend on the sequence $(B_m)$ at
all.)
\end{thm}

We will prove this theorem by using the Banach fixed point
theorem again. First we introduce the adequate Banach space.
To keep the notation as simple as possible, we define $\chi_n$ by
\begin{equation}                         \notag 
\chi_n(x)
= n^{-1/2} \varphi_{n \nu}(x)
+ n^{-d/2}\sum_{j=1}^{n/2} j\, \varphi_{j \nu}(x).
\end{equation}
In particular, we suppress the $\nu$-dependence of $\chi$. For
$\fAG\in\mathcal{S}$
we define the \ane{$\chi$-weighted} norm
\begin{equation*}
\normw{\fAG}:=\sup_{x\in\Zd}\abs{\AG_0(x)} + \sup_{n\in\N,x\in\Zd}
\frac{\abs{\AG_n(x)}}{\chi_n(x)},
\end{equation*}
whenever these suprema are finite. Finally we define the set
\begin{equation*}
\mathcal{W} \deff\left\{\fAG\in \mathcal{S}:
\normw{\fAG} <\infty \right\}.
\end{equation*}

Clearly $\mathcal{W}$ equipped with $\normw{.}$ is a Banach space.
We now have to determine the contraction operator for
the fixed point argument. This process will be more subtle than it was in
section \ref{chap-const}.

First we take an aperiodic probability measure $E$ in
$\mathcal{M}$ with covariance matrix
$\delta\cdot
\Id_d$. We can construct such a measure by taking the distribution
$\sist $ and shifting a small amount of the
probability to vectors of length $\stl+1$ (if $d \delta\ge \stv$) or
to zero (if $d \delta<\stv$). In this way, the range of $E$ is
bounded by $\step\deff \stl+1$.

We use this $E$ to define an appropriate contraction operator:
For
$\fAG\in\mathcal{S}$, the operator
$\fAG\mapsto\tifAG$ is given by the following recursion:
\begin{align*}
\tiAG_0 &= \AG_0, \qquad \text{and for }n\ge 1\\
\tiAG_n &= E*\tiAG_{n-1} + (\sist -E)*\AG_{n-1}
-\lambda\bigl[\rho \sist *\AG_{n-1} - \summ a_m\,B_m*\AG_{n-m}\bigr]
\end{align*}
for all sequences of signed measures $\fAG = (\AG_n)_{n\in\No}$ in
$\mathcal{S}$. Clearly $\tifAG$ is an element of $\mathcal{S}$.

The sequence $(A_n)$ that we have defined in \eqref{defA}, is obviously
a fixed point of $\sim$. We will show that it is asymptotically close to
the distribution
of a sum of i.i.d.\,random vectors with law $E$. We will use the following
lemmas:
\begin{lem}                                           \label{lem-schrittlok}
Let $\AG_n:= a_n E^{*n}$ and $\nu\ge 1/(2d)$ big enough. Then
\begin{equation*}
\normw{\tifAG-\fAG}\le \const  (1+\lambda \beta_\nu ),
\end{equation*}
where the
constant $\const $ depends on the dimension $d$ and $\stl$ only.
\end{lem}

\begin{lem}                                        \label{lem-kontraktlok}
Let $\fAG\in\mathcal{W}$ with $\AG_0=0$. For $\lambda$ small enough
there exists a parameter $\kappa\in(0,1)$ such that
\begin{equation*}
\normw{\tifAG}\le \kappa \normw{\fAG}.
\end{equation*}
\end{lem}

The proofs of the two lemmas are very similar. Straightforward calculation
allows us to rewrite $\tiAG_n$ in two different ways:
\begin{align}
\tiAG_n &= \AG_n - \sum_{l=1}^n
E^{*n-l} *\bigl[\AG_l-(1-\lambda\rho)\sist *\AG_{l-1} -            \label{D1}
\lambda\summl a_m\,B_m*\AG_{l-m}\bigr]\\
&\text{and}\notag\\
\tiAG_n &= E^{*n}*\AG_0 - \sum_{l=1}^n
\AG_{n-l}*\bigl[E^{*l}-(1-\lambda\rho) \sist *E^{*l-1}
-\lambda\summl a_m\,B_m*E^{l-m}\bigr].                         \label{D2}
\end{align}
These expressions can be viewed as perturbations of
the respective first term. The proofs of the lemmas now consist of error
estimates for the sums appearing as perturbation terms.

\begin{proof}[Proof of Lemma \ref{lem-schrittlok}]
According to \eqref{D1} it suffices to show that for all $n$
\begin{equation}\begin{split}                           \label{schritt1lok}
\sum_{l=1}^n\;&
\babs{a_l E^{*n}-(1-\lambda\rho)a_{l-1}\sist *E^{*n-1} -
\lambda\summl a_m\,B_m\,a_{l-m}*E^{*n-m}}\\
&\le  (1+\lambda \beta_\nu )\const \; \chi_n.
\end{split}\end{equation}
The strategy of the proof is to approximate the discrete distributions
$E^{*n}$ by fitting normal densities and use their Taylor expansion to obtain
the desired bounds. There are several error terms to control.

We use  Lemma \ref{lem-approx} to obtain an
approximation for $E^{*n}$. For $\nu'=\nu'(d,\step)$
large enough, the lemma yields
\begin{align}                                            \label{Endiff}
\babs{E^{*n}(x)-\bigl[1+n^{-1}P_4 (x/\sqrt{n}
)\bigr]\varphi_{n\delta}(x)} &\le \const  n^{-3/2}\varphi_{n \nu'}(x),
\end{align}
where $P_4$ is a polynomial of degree four. The coefficients of $P_4$ are
rational functions of the moments of $E$ up to order four. We now fix
$\nu$ as the maximum of $\sqrt{2}\nu'$ and $6/d$. In particular, together
with
\eqref{deltaest} this yields $\nu\ge 3\delta$, which we will use later
for the error estimates.

To simplify the notation, we use the following abbreviations:
\begin{align}                                       \notag 
\breve{\varphi}_n(x)&\deff
\bigl[1+n^{-1}P_4(x/\sqrt{n})\bigr]\varphi_{n\delta}(x)
\quad\text{denotes the function above},\\
Y_l &\deff a_l E-(1-\lambda\rho)a_{l-1}\sist   \qquad\text{and} \notag
\\X_l &\deff \babs{Y_l*\breve{\varphi}_{n-1}-
\lambda\summl[l\wedge n/2]
a_m\,a_{l-m}\,B_m*\breve{\varphi}_{n-m}}.  \label{Xl}
\end{align}
We split the left-hand side of \eqref{schritt1lok} into several parts, which
will be estimated separately:

\begin{subequations}                                 \label{schritt}
\begin{align}
\text{l.h.s. of \eqref{schritt1lok}}&\le
\sum_{l=1}^n\;X_l                                         \label{schritta}\\
&\quad+\sum_{l=1}^n\;\abs{Y_l}*\abs{\breve{\varphi}_{n-1} -
E^{*n-1}}\label{schrittb}\\
&\quad +\lambda \sum_{l=1}^n\;
\summl[{l\wedge n/2}]
a_m\,a_{l-m}\,\abs{B_m}*\abs{\breve{\varphi}_{n-m} - E^{*n-m}}
\label{schrittc}\\
&\quad + \lambda\sum_{l=n/2}^n\;
\sum_{m=n/2}^l a_m\,a_{l-m}\,\abs{B_m}*E^{*n-m}.
\label{schrittd}
\end{align}
\end{subequations}


Before we start to estimate the various sums, we want to state some
facts that we will use extensively in this proof.

Sums of the form $\sum_{l=1}^n l^{\varepsilon}$ for some real
$\varepsilon$ are normally estimated by majorizing them with the appropriate
integral $\int t^\varepsilon dt$. Double sums like $\sum_{l=1}^{n-1}
l^{\varepsilon} (n-l)^{\varepsilon'}$ are bounded by splitting them in
two parts $l\le n/2$ and $l\ge n/2$ and treating the halves separately.

Another  often used inequality yields an upper bound on the discrete folding
$\varphi_\eta*\varphi_\theta (x)\deff
\sum_{y\in\Zd}
\varphi_\eta(y) \varphi_\theta(x-y)$ of two normal densities. We have
\begin{equation*}
\varphi _{\eta}*\varphi _{\theta}(x) \le \const \varphi_{\eta+\theta}(x)
\end{equation*}
uniformly for all $x\in \Zd$ and  $\eta,\theta\ge 1/(2d)$. This inequality is
proven in Lemma
\ref{lem-faltung} in the appendix.

A last remark worth making is that for positive constants
$l$, $l'$ and $m$ with $l\le l'\le m\,l$ we have
\begin{equation*}
\varphi _{l \eta}(x) \le \const (m)\varphi_{l' \eta}(x)
\end{equation*}
for all $x\in \Zd$. This simple fact is obtained by bounding the first factor
and the exponential term in $\varphi_{l\eta}(x)$ separately.


Now we come back to, or rather start with, \eqref{schritt}. We go from bottom
to top and start with bounding \eqref{schrittd}.

A direct
consequence of Lemma \ref{lem-approx} is that $E^{*n} \le \const
\varphi_{n\nu'}$. So we obtain
\begin{align*}
\eqref{schrittd}
&=\lambda\sum_{l=n/2}^n\;
\sum_{m=n/2}^l a_m\,a_{l-m}\,\abs{B_m}*E^{*n-m}\\
&\le \lambda \beta_\nu \const  \sum_{l=n/2}^n\,\sum_{m=n/2}^l
\underbrace{m^{-d/2}}_{\le \const \, n^{-d/2}}
\sum_{k=1}^{m/2} k^{1-d/2}\underbrace{\varphi_{k
\nu}*\varphi_{(n-m)\nu'}}_ {\le \const  \varphi_{(k+n-m)\nu}}
\notag\\
&\le \lambda \beta_\nu \const \; n^{-d/2}\sum_{m=n/2}^n (n-m)
\sum_{k=1+n-m}^{n-m/2} (k-(n-m))^{1-d/2}\varphi_{k \nu}
\notag\\
&\le \lambda \beta_\nu \const \; n^{-d/2} \sum_{k=1}^{3n/4}\varphi_{k \nu}
\sum_{m=n/2\vee n-k+1}^n (k-(n-m))^{1-d/2}\underbrace{(n-m)}_{\le k}
\notag\\
&\le \lambda \beta_\nu \const \; \chi_n.
\end{align*}
In the last step we used the fact, that --- in five and more dimensions
--- the sum $n^{-d/2}\sum_{k=1}^{n} k \varphi_{k \nu}$ is bounded above by
$\chi_n$, which can easily be seen by splitting the sum.
We use \eqref{Endiff} to see
\begin{align*}
\eqref{schrittc}
&=\lambda \sum_{l=1}^n\;
\summl[{l\wedge n/2}]
a_m\,a_{l-m}\,\abs{B_m}*\abs{\breve{\varphi}_{n-m} - E^{*n-m}}\notag\\
&\le \lambda \const  \sum_{l=1}^n\;
\summl[{l\wedge n/2}]\,(n-m)^{-3/2} \abs{B_m}*\varphi_{(n-m) \nu'}\\
&\le \lambda \beta_\nu \const  n^{-1/2}\summl[n/2]  m^{-d/2}
\sum_{k=1}^{m/2} k^{1-d/2} \underbrace{\varphi_{k \nu} *\varphi_{(n-m)
\nu'}} _{\le \const  \varphi_{(n-(m-k)) \nu}}\\
&\le \lambda \beta_\nu \const  n^{-1/2}\varphi_{n \nu},
\end{align*}
where we used in the last line that $n/2 \le n-(m-k)\le n$ (and therefore
$\varphi_{(n-(m-k)) \nu}\le \const  \varphi_{n \nu}$).
Furthermore we have
\begin{align*}
\eqref{schrittb}
&=\sum_{l=1}^n\;\abs{Y_l}*\abs{\breve{\varphi}_{n-1} -
E^{*n-1}}\\
&\le \const (n-1)^{-3/2}
\sum_{l=1}^n\;\underbrace{\abs{Y_l}*\varphi_{(n-1)
\nu'}}_{\le K \varphi_{\sqrt{2}(n-1)\nu'}}
\le \const \,n^{-1/2}\varphi_{n \nu},
\end{align*}
since $Y_l= a_l E-(1-\lambda\rho)a_{l-1}\sist$ is a signed measure of bounded
steplength.

Now we come to \eqref{schritta}. Recall the definition of
$X_l$ in \eqref{Xl}, that is
\begin{equation*}
X_l \deff \babs{Y_l*\breve{\varphi}_{n-1}-
\lambda\summl[l\wedge n/2]
a_m\,a_{l-m}\,B_m*\breve{\varphi}_{n-m}}.
\end{equation*}
We analyze the terms in $X_l$ separately, using
Lemma
\ref{lem-A} and Lemma \ref{lem-B}. We write $P_2(z)$ for the polynomial
$z^2/\delta -d$.

        From \eqref{A} we obtain
\begin{align}                                           \label{teil1}
Y_l * \breve{\varphi}_{n-1}(x)
&= \bigl[a_l-(1-\lambda\rho)a_{l-1}\bigr]\cdot\breve{\varphi}_{n-1}(x) \notag\\
            &\quad+ \Bigl[\frac{a_l
d\delta-(1-\lambda\rho)a_{l-1}\stv}{2(n-1)d\delta}\Bigr]
\cdot P_2(x/\sqrt{n-1})\varphi_{(n-1)\delta}(x) \notag\\
&\quad+ R_4^{Y_l}(n-1,x)
                \; + (n-1)^{-1}R_2^{Y_l}(n-1,x;P_4).
\end{align}

Now we apply first \eqref{B} and then \eqref{A} to obtain
\footnote{In the following formulas the dots denote the argument with respect
to which we fold.}
\begin{align}                                           \label{teil2}
B_m*\breve{\varphi}_{n-m}(x)
&=B_m*\breve{\varphi}_{n-1}(x) -\frac{m-1}{2(n-1)}
\bigl[B_m*P_2(./\sqrt{n-1})\varphi_{(n-1)\delta}\bigr](x)\notag\\
&\quad +
\bigl[B_m*S_2^{m-1}(n-1,.)\bigr](x)\notag\\
&\quad +
\bigl[B_m*S_1^{m-1}(n-1,.;(n-1)^{-1}P_4)\bigr](x)\notag\\[3mm]
&= b_m \cdot \breve{\varphi}_{n-1}(x)\notag\\
&\quad
+\Bigl[\frac{\ulb_m}{2(n-1)d \delta} - \frac{(m-1)b_m}{2(n-1)}\Bigr]\cdot
P_2(x/\sqrt{n-1})\varphi_{(n-1)\delta}(x) \notag\\
&\quad +R_4^{B_m}(n-1,x)  +(n-1)^{-1}R_2^{B_m}(n-1,x;P_4)\notag\\
&\quad -\frac{m-1}{2(n-1)}R_2^{B_m}(n-1,x;P_2)
           + \bigl[B_m*S_2^{m-1}(n-1,.)\bigr](x) \notag\\ &\quad +
\bigl[B_m*S_1^{m-1}(n-1,.;(n-1)^{-1}P_4)\bigr](x).
\end{align}

We insert \eqref{teil1} and \eqref{teil2} into \eqref{Xl} and obtain
\begin{align}                                           \label{schritt3lok}
X_l(x) \le
I_l\cdot \underbrace{\breve{\varphi}_{n-1}(x)}_{\le \const
\,\varphi_{2n\delta}(x)} +
J_l\cdot\underbrace{\abs{P_2(x/\sqrt{n-1})}\varphi_{(n-1)\delta}} _{\le
\const \,\varphi_{2n\delta}(x)} + \,R_l(x),
\end{align}
with
\begin{align*}
I_l &=\babs{a_l-(1-\lambda\rho)a_{l-1} -\lambda\summl[l\wedge n/2] a_m b_m
a_{l-m}}\qquad\qquad \text{and}\\
J_l &=\frac{1}{2(n-1)d \delta}\Babs{a_ld \delta-(1-\lambda\rho)a_{l-1}\stv
-\lambda\summl[l\wedge n/2] a_m a_{l-m} (\ulb_m  -d
\delta(m-1)b_m)}.
\end{align*}

Using the recursion formula for $a_l$  and the  decay rate
of $b_m=O(m^{-d/2})$ (see \eqref{bmdecay}), we obtain
\begin{align}                                               \label{Illok}
I_l &
\le \lambda\sum_{m=n/2}^l a_m \abs{b_m} a_{l-m}\le \lambda \beta_\nu
\const  n^{-3/2}.
\end{align}

Recall also (see Corollary \ref{cor:speed}) that $b_m=O(m^{-5/2})$ implies
$(\Delta a)_n=O(n^{-3/2})$. More precisely, we have
$\abs{\alpha-a_n}\le \lambda
\beta_\nu \const  n^{-1/2}$ and therefore
\begin{align}
J_l &=\frac{1}{2(n-1)d \delta}\Babs{a_ld \delta-(1-\lambda\rho)a_{l-1}\stv
-\lambda\sum_{m=1}^{l\wedge n/2} a_m a_{l-m} (\ulb_m  -d
\delta(m-1)b_m)}\notag\\
&\le \frac{\alpha}{2(n-1)d \delta}\Babs{d
\delta-(1-\lambda\rho)\stv -\lambda\sum_{m=1}^{l\wedge n/2} a_m\ulb_m
+d \delta\lambda\sum_{m=1}^{l\wedge n/2}
(m-1)a_m b_m} \notag\\
&\qquad + \lambda \beta_\nu \const  n^{-1} l^{-1/2}\notag\\
&\le\frac{\alpha}{2(n-1)d \delta}\babs{\underbrace{d \delta-(1-\lambda\rho)
\stv
-\lambda\tau  +d \delta\lambda\sigma}_{=0}}+ \lambda \beta_\nu \const  n^{-1}
l^{-1/2}\notag\\ &\le \lambda \beta_\nu \const  n^{-1}
l^{-1/2},
\label{Jllok}
\end{align}
where we used \eqref{bmdecay} again and the fact that $\nu=\nu(d,\step)$ has
already been fixed and can thus be bounded by a constant $\const $.
The error term $R_l$ in \eqref{schritt3lok} is given by
\begin{subequations}\label{fehl-lok}
\begin{align}
R_l(x)&= \abs{R_4^{Y_l}(n-1,x)}
                + (n-1)^{-1}\abs{R_2^{Y_l}(n-1,x;P_4)}\label{fehl-alok}\\
&\quad +\lambda\summl[l\wedge n/2] a_m a_{l-m}
\abs{R_4^{B_m}(n-1,x)}\label{fehl-blok}\\ &\quad
+\lambda(n-1)^{-1}\summl[l\wedge n/2] a_m a_{l-m}
\abs{R_2^{B_m}(n-1,x;P_4)}\label{fehl-clok}\\
&\quad +\lambda(2(n-1))^{-1}\summl[l\wedge n/2] (m-1) a_m
a_{l-m}\abs{R_2^{B_m}(n-1,x;P_2)}\label{fehl-dlok}\\
&\quad + \lambda\summl[l\wedge n/2] a_m
a_{l-m} \abs{S_2^{m-1}(n-1,.)*B_m(x)}\label{fehl-elok}\\
&\quad + \lambda\summl[l\wedge n/2]
a_m a_{l-m} \abs{S_1^{m-1}(n-1,.;(n-1)^{-1}P_4)*B_m(x)}. \label{fehl-flok}
\end{align}
\end{subequations}
We use the local error estimates in lemmas \ref{lem-A} and \ref{lem-B} to
estimate the various terms. We have
\begin{align*}
\abs{R_4^{Y_l}(n-1,x)}&\le \const  n^{-2}\varphi_{2n\delta}(x)
\qquad\text{by \eqref{R},}
\end{align*}
since $\abs{Y_l}$ has bounded steplength.
Analogously we can show the same decay for the second term of
\eqref{fehl-alok}.

On the other hand, we bound \eqref{fehl-blok} using \eqref{R'} to obtain
\begin{align*}
\lambda \const &\summl[l\wedge n/2]\babs{R_4^{B_m}(n-1,x)}\notag\\
&\le \lambda \const \summl[l\wedge n/2]\; n^{-2}\int_0^1 \!\!ds
\,\,\sum_{z\in\Zd} z^4\,
\abs{B_m(z)}\;\varphi_{\sqrt{2}(n-1)\delta}(x-sz)\notag\\
&\le \lambda \beta_\nu
\const  n^{-2} \summl[l\wedge n/2] m^{-d/2}\sum_{k=1}^{m/2}
k^{1-d/2}\\
&\hspace{3cm}\times\int_0^1 \!\!\!ds  \sum_{z\in\Zd}
\underbrace{z^4\;\varphi_{k \nu}(z)}_{\le \const  k^2\varphi_{2{k \nu}}(z)}\;
\underbrace{\varphi_{\sqrt{2}(n-1)\delta}(x-sz)}_{\le
\const s^{-d}\,\varphi_{\sqrt{2}(n-1)\delta/s^2}(x/s-z)}\notag\\
&\le \lambda \beta_\nu \const  n^{-2} \summl[l\wedge n/2]
m^{-d/2}\sum_{k=1}^{m/2}
k^{3-d/2} \int_0^1 \!\!ds\;
\underbrace{s^{-d}\varphi_{2{k \nu}+\sqrt{2}(n-1)\delta/s^2}(x/s)}%
_{\le \const \varphi_{2{k \nu}s^2+\sqrt{2}(n-1)\delta}(x)}\notag\\
&\le \lambda \beta_\nu \const  n^{-2} \varphi_{n
\nu}(x)\underbrace{\summl[l\wedge n/2]
m^{-d/2}\sum_{k=1}^{m/2} k^{3-d/2}}_{\le \const  \,\log(n)}
\notag \\ &\le \lambda \beta_\nu \const  n^{-3/2}\varphi_{n \nu}(x),
\end{align*}
where we used $\nu \ge 3\delta$ in the second to last step.
Analogously we can show that \eqref{fehl-clok} and \eqref{fehl-dlok} are
bounded by $\lambda \beta_\nu \const  n^{-2}\varphi_{n \nu}(x)$ and
$\lambda \beta_\nu \const  n^{-3/2}\varphi_{n \nu}(x)$, respectively.
We bound \eqref{fehl-elok} with \eqref{S} by
\begin{align*}
\lambda \beta_\nu  \const  &\summl[l\wedge n/2] m^{2-d/2}
\underbrace{(n-m)^{-2}
\bigl(\frac{n}{n-m}\bigr)^{d/2}}_{\le \const  n^{-2}}
\sum_{k=1}^{m/2}k^{1-d/2}\varphi_{\sqrt{2}(n-1)\delta}*\varphi_{k \nu} (x)\\
&\le \lambda \beta_\nu \const  n^{-2}\summl[l\wedge n/2] m^{2-d/2}
\sum_{k=1}^{m/2}k^{1-d/2}\varphi_{{k \nu}+\sqrt{2}(n-1)\delta} (x)\\
&\le \lambda \beta_\nu \const  n^{-2}\varphi_{n \nu}(x)\summl[l\wedge n/2]
m^{2-d/2}
\underbrace{\sum_{k=1}^{m/2}k^{1-d/2}}_{\le \const }\\
&\le \lambda \beta_\nu \const  n^{-2}l^{1/2}\varphi_{n \nu}(x)\le \lambda
\beta_\nu \const  n^{-3/2}\varphi_{n \nu}(x),
\end{align*}
where we used $\nu\ge 3\delta$ to obtain the second to last line.
For \eqref{fehl-flok} an even better bound without the
$l^{1/2}$-term results analogously.

Combining the different estimates, we can bound $R_l$ by
\begin{align*}
R_l
&\le (1+\lambda \beta_\nu) \const \, n^{-3/2}\varphi_{n \nu}.
\end{align*}
Considering this together with \eqref{Illok} and \eqref{Jllok},
we can finally bound \eqref{schritta}:
\begin{align*}
\eqref{schritta}
&=\sum_{l=1}^n\;X_l \\
&\le
\sum_{l=1}^n\,
I_l \;\varphi_{2n\delta}(x)+ J_l\;\varphi_{2n\delta}(x) + \,R_l(x)\\
&\le  (1+\lambda \beta_\nu ) \const  \,\varphi_{n
\nu}(x)\sum_{l=1}^n\Bigl[n^{-3/2} + n^{-1} l^{-1/2}\Bigr]\\
&\le (1+\lambda \beta_\nu ) \const  \,n^{-1/2} \varphi_{n \nu}(x).
\end{align*}
This proves Lemma \ref{lem-schrittlok}.
\end{proof}

\begin{proof}[Proof of Lemma \ref{lem-kontraktlok}]
Let $\fAG$ with $\AG_0=0$ and $\normw{\fAG}:=
\sup_{n,x} \chi_n(x)^{-1} \abs{\AG_n(x)}< \infty$ be given.

We want to prove that for small enough $\lambda$ we have
$\normw{\tifAG}\le\kappa\normw{\fAG}$ for some $\kappa<1$.
According to \eqref{D2}, it is sufficient to show that for all natural $n$
\begin{align*}
\sum_{l=1}^{n-1}
\babs{\AG_{n-l}}*\babs{E^{*l}-(1-\lambda\rho) \sist*E^{*l-1}
-\lambda\summl a_m\,&B_m*E^{l-m}}\\&\le \lambda (1+\beta_\nu) \const
\chi_n\norm{\fAG}.
\end{align*}
Once we know this, we can fix $\const $ and choose $\lambda$ so small that
$\lambda (1+\beta_\nu) \const $ is strictly smaller than one.
Since $\babs{\AG_{n-l}}\le \norm{\fAG}\cdot\chi_{n-l}$, it suffices to
show that
\begin{align}                                         \label{kontrakt1lok}
\sum_{l=1}^{n-1} \chi_{n-l}*\babs{\bigl(E-(1-\lambda\rho)S\bigr)
*E^{*l-1} -\lambda\summl a_m\,B_m*E^{*l-m}}&\le
\lambda(1+\beta_\nu) \const  \chi_n.
\end{align}
The proof of this estimate is somewhat tedious, but  similar to the
preceding one.  We take the same $\nu$ and we will split the sum in
the very same way as before. Again we
use some abbreviations to
         keep the notation as readable as possible:
\begin{align*}
\breve{\varphi}_n(x)&\deff
\bigl[1+n^{-1}P_4(x/\sqrt{n})\bigr]\varphi_{n\delta}(x)
\quad\text{as before},\notag\\
Y &\deff E-(1-\lambda\rho) \sist  \qquad\text{and}      
\\
X'_l &\deff \babs{Y*\breve{\varphi}_{l-1}- \lambda
\summl[l/2] a_m\,B_m*\breve{\varphi}_{l-m}}.  
\end{align*}
The left hand side of \eqref{kontrakt1lok} is split in the following parts,
which will be estimated separately:
\begin{subequations}
\begin{align}
\text{l.h.s. of \eqref{kontrakt1lok}}&\le
\sum_{l=1}^{n-1}\;\chi_{n-l} *X'_l
\label{kontrakta}\\
&\quad+\sum_{l=1}^n\;\chi_{n-l}*\abs{Y}*\abs{\breve{\varphi}_{l-1} -
E^{*l-1}}                                                 \label{kontraktb}\\
&\quad +\lambda\sum_{l=1}^{n-1}\;\chi_{n-l}*
\summl[l/2] a_m\;\abs{B_m}*\abs{\breve{\varphi}_{l-m} - E^{*l-m}}
\label{kontraktc}\\
&\quad +\lambda \sum_{l=1}^{n-1}\;\chi_{n-l}*
\sum_{m=l/2}^l a_m\;\abs{B_m}*E^{*l-m}.        \label{kontraktd}
\end{align}
\end{subequations}

We again start with the last term. First we split $\chi_{n-l}$ into
$(n-l)^{-1/2} \varphi_{(n-l)\nu}$
and $(n-l)^{-d/2}\sum_{j=1}^{(n-l)/2}j\,\varphi_{j \nu}$. We treat the
resulting parts of \eqref{kontraktd} separately. For the first part this
leads to
\begin{align}                               \notag 
\sum_{l=1}^{n-1} & (n-l)^{-1/2} \varphi_{(n-l)\nu}*
\sum_{m=l/2}^{l}a_m\,\abs{B_m}*E^{*l-m}\notag\\
&\le \beta_\nu \const  \sum_{l=1}^{n-1}  (n-l)^{-1/2}
\sum_{m=l/2}^{l}m^{-d/2}\sum_{k=1}^{m/2}
k^{1-d/2}\,\underbrace{ \varphi_{(n-l)\nu}*\varphi_{k
\nu}*\varphi_{(l-m)\nu'}}_{\le \const
\varphi_{(n-m+k)\nu}}\notag\\
&\le \beta_\nu \const  n^{-1/2} \varphi_{n \nu} \underbrace{\sum_{l=1}^{n/2}
\sum_{m=l/2}^{l}m^{-d/2}}_{\le \const }\underbrace{\sum_{k=1}^{m/2}
k^{1-d/2}}_{\le \const }\notag\\
&\quad + \beta_\nu \const  n^{-d/2}\sum_{l=n/2}^{n-1}  (n-l)^{-1/2}
\sum_{m=l/2}^{l}\;\sum_{k=1}^{m/2}
k^{1-d/2} \varphi_{(n-m+k) \nu}\notag\\
&\le \beta_\nu \const  \Bigl[ n^{-1/2} \varphi_{n \nu}
           + n^{-d/2}\!\!\!
\sum_{m=n/4}^{n-1}\;\sum_{k=n-m+1}^{n-m/2}
\!\!\!\!(k-(n-m))^{1-d/2} \varphi_{k \nu}
\underbrace{\sum_{l=m}^{n-1} (n-l)^{-1/2}}_{\le n-m\,\,\le\, k}\Bigr]
\notag\\
&\le \beta_\nu \const  \Bigl[ n^{-1/2} \varphi_{n \nu}
           + n^{-d/2}
\sum_{k=1}^{7n/8}k\, \varphi_{k \nu}
\underbrace{\sum_{m=n-k+1}^{n-1}(k-(n-m))^{1-d/2}}_{\le \const }
\Bigr]
\notag\\
&\le \beta_\nu \const  \chi_n.     \notag %
\end{align}
For the second part we split the sum and find
\begin{align}                                \notag 
\sum_{l=1}^{n/2}  &(n-l)^{-d/2} \sum_{j=1}^{(n-l)/2} j\, \varphi_{j \nu}
*\sum_{m=l/2}^{l}a_m\,\abs{B_m}*E^{*l-m}\notag\\
& \le \beta_\nu \const  \, n^{-d/2}\sum_{l=1}^{n/2} \, \sum_{j=1}^{(n-l)/2}j\,
\sum_{m=l/2}^{l} m^{-d/2} \sum_{k=1}^{m/2} k^{1-d/2}\varphi_{(j+k+l-m)\nu}
\notag\\
& \le \beta_\nu \const  \, n^{-d/2} \sum_{m=1}^{n/2} m^{-d/2}
\sum_{l=m}^{n/2\wedge 2m}
\, \sum_{k=1}^{m/2} k^{1-d/2}\sum_{j=1+k}^{(n-l)/2+k}
j\,\varphi_{(j+l-m)\nu} \notag\\
& \le \beta_\nu \const  \, n^{-d/2} \sum_{m=1}^{n/2} m^{-d/2}
\sum_{l=m}^{n/2\wedge 2m}
\,\sum_{j=1+l-m}^{(n+l-m)/2}
j\,\varphi_{j \nu} \, \sum_{k=1}^{m/2} k^{1-d/2}
\notag\\
& \le \beta_\nu \const  \, n^{-d/2} \sum_{j=1}^{3n/4}
j\,\varphi_{j \nu} \,  \sum_{m=1}^{n/2} m^{1-d/2}
\notag\\
& \le \beta_\nu \const  \chi_n, \notag %
\end{align}
and finally
\begin{align}                                   \notag 
\sum_{l=n/2}^{n-1}  &(n-l)^{-d/2} \sum_{j=1}^{(n-l)/2} \underbrace{j}_{\le
n-l}\,
\varphi_{j \nu} *\sum_{m=l/2}^{l}a_m\,\abs{B_m}*E^{*l-m}\notag\\
& \le \beta_\nu \const  \,\sum_{l=n/2}^{n-1} (n-l)^{1-d/2}\,
\sum_{m=l/2}^{l}\; m^{-d/2}\sum_{k=1}^{m/2} k^{1-d/2}\sum_{j=1}^{(n-l)/2}
\,\varphi_{(j+k+l-m)\nu}
\notag\\
& \le \beta_\nu \const  \,n^{-d/2}\sum_{l=n/2}^{n-1} (n-l)^{1-d/2}\,
\sum_{m=0}^{l/2}\; \sum_{k=1}^{(l-m)/2} k^{1-d/2}\sum_{j=1}^{(n-l)/2}
\,\varphi_{(j+k+m)\nu}
\notag\\
& \le \beta_\nu \const  \,n^{-d/2}\sum_{l=n/2}^{n-1} (n-l)^{1-d/2}\,
\sum_{m=0}^{l/2}\;\sum_{j=1+m}^{(n+m)/2}
\,\varphi_{j \nu}
\notag\\
& \le \beta_\nu \const  \,n^{-d/2}\sum_{l=n/2}^{n-1} (n-l)^{1-d/2}\;
\sum_{j=1}^{(2n+l)/4}\,\varphi_{j \nu}\sum_{m=0}^{j-1}
\notag\\
& \le \beta_\nu \const  \,n^{-d/2}\sum_{j=1}^{3n/4}j\,\varphi_{j \nu}
\notag\\
& \le \beta_\nu \const  \chi_n. \notag %
\end{align}
Thus $\eqref{kontraktd}\le \lambda\beta_\nu \const  \chi_n$ is shown.
In order to estimate \eqref{kontraktc}, we use
\eqref{Endiff} again, obtaining
\begin{align*}
\summl[l/2]&
a_m\;\abs{B_m}*\abs{\breve{\varphi}_{l-m} - E^{*l-m}}\notag\\
&\le \beta_\nu \const  \summl[l/2](l-m)^{-3/2} m^{-d/2}
\sum_{k=1}^{m/2} k^{1-d/2} \underbrace{\varphi_{k \nu}*\varphi_{(l-m)
\nu'}}_{\le \const  \varphi_{l \nu}} \notag\\
&\le \beta_\nu \const  l^{-3/2}\varphi_{l \nu},
\end{align*}
and therefore
\begin{align}                                         \label{kontrakt8lok}
\eqref{kontraktc}
&= \lambda\sum_{l=1}^{n-1}\;\chi_{n-l}*
\summl[l/2] a_m\;\abs{B_m}*\abs{\breve{\varphi}_{l-m} - E^{*l-m}}\notag\\
&\le \lambda\beta_\nu \const  \sum_{l=1}^{n-1}\;l^{-3/2}\,\chi_{n-l}
*\varphi_{l\nu}\notag\\
&\le \lambda\beta_\nu \const  \sum_{l=1}^{n-1}\;l^{-3/2}\,
\Bigl[ (n-l)^{-1/2}\varphi_{n \nu}+ (n-l)^{-d/2}\sum_{j=1}^{(n-l)/2}
j\, \varphi_{(j+l)\nu}\Bigr]\notag\\
&\le \lambda\beta_\nu \const  n^{-1/2}\varphi_{n \nu}
+ \lambda\beta_\nu \const  n^{-d/2}\sum_{l=1}^{n/2}\;l^{-3/2}\,
\sum_{j=1}^{(n-l)/2} j\, \varphi_{(j+l)\nu}\notag\\
&\qquad\qquad\qquad\qquad\,+\lambda\beta_\nu \const
\sum_{l=n/2}^{n-1}\;l^{-3/2}\,(n-l)^{-d/2}\sum_{j=1}^{(n-l)/2} j\,
\underbrace{\varphi_{(j+l)\nu}}_{\le \const  \varphi_{n \nu}}
\notag\\
&\le \lambda\beta_\nu \const  \chi_n.
\end{align}

Analogously, we have
\begin{align*}
\eqref{kontraktb}&=
\sum_{l=1}^n\;\chi_{n-l}*\abs{Y}*\abs{\breve{\varphi}_{l-1} -
E^{*l-1}}\le \lambda\beta_\nu \const  \chi_n.
\end{align*}

The remaining part, \eqref{kontrakta}, will be treated like
\eqref{schritta} in the proof of the last lemma, using Lemma
\ref{lem-A} and Lemma
\ref{lem-B}. Again we write $P_2(z)$ for the polynomial
$z^2/\delta -d$. From \eqref{A} we have
\begin{align}                            \notag 
Y*\breve{\varphi}_{l-1}(x)
&= \lambda\rho\cdot\breve{\varphi}_{l-1}(x) \notag\\
            &\quad+ 
\Bigl[\frac{d\delta-(1-\lambda\rho)\stv}{2(l-1)d\delta}\Bigr]
\cdot P_2(x/\sqrt{l-1})\varphi_{(l-1)\delta}(x) \notag\\
&\quad+ R_4^{Y}(l-1,x)
                \; + (l-1)^{-1}R_2^{Y}(l-1,x;P_4). \notag %
\end{align}

For the second sum, we use \eqref{teil2} after having replaced $n$ by $l$
everywhere. We obtain
\begin{align}                                           \label{kontrakt6lok}
X'_l(x) \le
I'_l\cdot \underbrace{\breve{\varphi}_{l-1}(x)}_{\le \const
\,\varphi_{2l\delta}(x)} +
J'_l\cdot\underbrace{\abs{P_2(x/\sqrt{l-1})}\varphi_{(l-1)\delta}} _{\le
\const \,\varphi_{2l\delta}(x)} + \,R'_l(x)
\end{align}
with
\begin{align*}
I'_l &=\babs{\lambda\rho -\lambda\summl[l/2] a_m b_m}\le \lambda\beta_\nu
\const  l^{-3/2}
\qquad\qquad \text{and}\\
J'_l &=\frac{1}{2(l-1)d \delta}\Babs{d\delta+(1-\lambda\rho)\stv
-\lambda\summl[l/2] a_m (\ulb_m -d \delta(m-1)b_m)}
\le \lambda\beta_\nu \const  l^{-3/2},
\end{align*}
which results from an argument very similar to that which led to
\eqref{Illok} and \eqref{Jllok}.
This time the error term $R'_l$ in \eqref{kontrakt6lok} is given by
\begin{subequations}
\begin{align}
R'_l(x)&= \abs{R_4^{Y}(l-1,x)}
                + (l-1)^{-1}\abs{R_2^{Y}(l-1,x;P_4)}\label{fehl-aalok}\\
&\quad +\lambda\summl[l/2] a_m
\abs{R_4^{B_m}(l-1,x)}\label{fehl-bblok}\\ &\quad
+(l-1)^{-1}\lambda\summl[l/2] a_m
\abs{R_2^{B_m}(l-1,x;P_4)}\label{fehl-cclok}\\
&\quad +(l-1)^{-1}\lambda\summl[l/2] (m-1) a_m
\abs{R_2^{B_m}(l-1,x;P_2)}\label{fehl-ddlok}\\
&\quad + \lambda\summl[l/2] a_m
           \abs{S_2^{m-1}(l-1,.)*B_m(x)}\label{fehl-eelok}\\
&\quad +\lambda\summl[l/2]
a_m a_{l-m} \abs{S_1^{m-1}(l-1,.;(l-1)^{-1}P_4)*B_m(x)}. \label{fehl-fflok}
\end{align}
\end{subequations}
These terms are bounded in exactly the same way as the corresponding ones in
formula \eqref{fehl-lok}. So we obtain
\begin{align*}
R'_l
&\le \lambda (1+\beta_\nu ) \const \, l^{-3/2}\varphi_{l \nu}.
\end{align*}
Combining these estimates we get
\begin{align*}
\eqref{kontrakta}
&\le\lambda(1+\beta_\nu ) \const \, \sum_{l=1}^{n-1}
l^{-3/2}\;\chi_{n-l}*\varphi_{l
\nu}
\le \lambda(1+\beta_\nu )\const  \chi_n
\end{align*}
as in \eqref{kontrakt8lok}.
\end{proof}

Now we have the necessary tools to prove Theorem \ref{thm-local}:

\begin{proof}[Proof of Theorem \ref{thm-local}]
First recall the space
$\mathcal{W}= \{\fAG\in\mathcal{S}: \normw{\fAG}<\infty\}$ and write
$\mathcal{W}_0$ for $\{\fAG\in\mathcal{W}: G_0\equiv 0\}$. For the
moment, we refer to the sequence
$(a_n E^{*n})_{n\in\No}$ as ${\bf{E}}$. The following considerations always
assume that $\lambda$ is small enough. Lemma
\ref{lem-schrittlok} yields
$\bigl({\bf{E}}-\widetilde{\bf{E}}\bigr)\in\mathcal{W}_0$. From Lemma
\ref{lem-kontraktlok} we know that $\sim$ is a contraction on
$\mathcal{W}_0$. Since $\sim$ is linear,
the Banach fixed point theorem now yields
the existence of a unique fixed point in ${\bf{E}}+\mathcal{W}_0$. With
other words, we  have a unique
sequence $\fA$ of symmetric measures on $\Zd$ with
\begin{enumerate}
\item
$A_0=\delta_0$,
\item
$\tiaa_n(x)=A_n(x)$\quad for all $n\in\N$, $x\in\Zd$,
\item
$\normw{{\bf{E}}-\fA}\le \const  (1+\lambda\beta_\nu)$.
\end{enumerate}
This sequence obviously {\em is} the sequence $(A_n)_{n\in\No}$
defined in \eqref{defA}.
Since as a direct consequence of Lemma \ref{lem-approx} we have
\begin{equation*}
\abs{a_n E^{*n}(x)- a_n \varphi_{n\delta}(x)} \le \const  n^{-1/2}
\varphi_{n\nu} (x),
\end{equation*}
the estimate \eqref{local} now results by choosing
$\lambda_0=\lambda_0(d,\beta_\nu)$ small enough.
\end{proof}

To end this section, we briefly discuss the periodic case.
Assume that $S$ and the sequence $(B_m)_{m\in\No}$ are {\em
two-periodic} (that is,
$B_m(x)=0$ whenever $m$ and $\norme{x}$ do not have the same parity). By
\eqref{defA}, the periodicity transfers to the whole sequence $(A_n)$. In
order to use the same arguments as in the aperiodic case, we have to define
a periodic probability measure $E'$ so that we can approximate $A_n$ by
$(E')^{*n}$. If the diffusion constant $\delta$ (see \eqref{deltadef}) is
greater or equal to $1/d$, this can be done easily.

So let $S\in\mathcal{M}$ be a two-periodic distribution  of
bounded range less or equal $\stl$. Furthermore let $(B_m)\in\mathcal{S}$ be
a two-periodic sequence which obeys
\eqref{Bmlokprop} and has parameter $\delta$ (resulting from
this sequence by \eqref{deltadef}) greater or equal to $1/d$. Then we have the
following theorem.
\begin{thm}[Local Estimates, two-periodic case]        \label{thm-localper}
Under the above assumptions, the sequence $(A_n)$
defined by
\eqref{defA}, has the following property: There exist $\lambda_0>0$ small
enough and $\nu>0$ big enough such that for all $\lambda\in(0,\lambda_0)$
\begin{align*}
\abs{A_n(x) - 2a_n\varphi_{n\delta}(x)} &\le \const \Bigl[
n^{-1/2}\varphi_{n \nu}(x) + n^{-d/2}\sum_{j=1}^{n/2} j\,\varphi_{j
\nu}(x)\Bigr],
\end{align*}
where $n$ is taken to have the same parity as $\norme{x}$ and $\const$ is a
positive constant depending only on
$d$ and $\stl$.
\end{thm}

\begin{proof}
Replace the aperiodic $E$ by the periodic $E'$ everywhere. Instead of using
\eqref{Endiff}, apply a periodic version of Lemma \ref{lem-approx}, namely
\begin{align*}
\babs{E'{}^{*n}(x)-2\bigl[1+n^{-1}P_4 (x/\sqrt{n}
)\bigr]\varphi_{n\delta}(x)} &\le \const  n^{-3/2}\varphi_{n \nu}(x),
\end{align*}
whenever $\norme{x}$ and $n$ have the same parity. The rest of the
proof for the
aperiodic case  carries over word by word (with the constants
suitably adapted).
\end{proof}

Note that if $\delta$ is smaller than
$1/d$, we have a problem with our construction. A symmetric and
rotationally invariant, twoperiodic probability measure on the
lattice $\Zd$ with
variance smaller than $1/d$ simply does not exist. Possibly this case can be
covered by choosing a more delicate contraction operator.

The weakly self-avoiding walk is spreading faster rather
than slower compared with the simple random walk. Therefore the
speed of its diffusion is greater or equal to the one of the simple random
walk. This means
$\delta\ge 1/d$. Instead of giving this heuristic argument one can calculate
for small
$\lambda$ the leading term in $\delta$ (use \eqref{deltaid} and \eqref{Pim}
as well as
\eqref{Pim1}), which is larger than $1/d$.

Therefore it will be possible to apply Theorem \ref{thm-localper} to the
weakly
self-avoiding walk as soon as we have shown that the lace functions have the
desired decay property \eqref{Bmlokprop}. This is the content of the next
section.
\section{Application to the Weakly Self-Avoiding Walk} \label{chap-saw}

We now come back to the specific context of the weakly self-avoiding walk,
where the main objects of study are the two-point functions
$C_n$ with total mass $c_n$. Recall that they satisfy the lace expansion
formula \eqref{recursion}, that is
\begin{equation*}
         C_n = 2d(D*C_{n-1}) + \sum_{m=2}^n \Pi_m*C_{n-m}.
\end{equation*}

\subsection{Decay Behavior}

We still have to show that the methods in the last chapters can in fact be
applied to the weakly self-avoiding walk. This means that we have to show
that $\Pi_m/(\lambda c_m)$ really has the behavior we assumed for $B_m$ in
section
\ref{chap-local}.

\begin{lem}                                                \label{lem-behave}
There are positive constants $\lambda_0$, $\nu$ and $\decc$, such that for all
$\lambda\in(0, \lambda_0)$ we have
\begin{align}                                  \notag 
\frac{\abs{\Pi_m(x)}}{\lambda c_m}&\le
\decc m^{-d/2}\sum_{k=1}^{m/2}k^{1-d/2}\,\varphi_{k\nu}(x).
\end{align}
\end{lem}

\begin{proof}
The argument uses induction on $m$. As in section \ref{chap-def}, we define
$B_m(x)\deff
\Pi_m(x)/(\lambda c_m)$. We freely use the results of appendix $\ref{chap-le}$.

We have $\Pi_1(x)=0$ and hence $B_1\equiv 0$ too.
Now consider $\Pi_2$: There is only one lace of length two, namely
$\{02\}$. Therefore we have by \eqref{Pim} and \eqref{PimNdef}
\begin{align}                                    \notag 
\Pi_2(x) &= -\lambda \sumsx[=2] U_{02}(\omega) = -2d\lambda \delta_{0x}.
\end{align}
On the other hand, we calculate easily (recall \eqref{Cndef} and
\eqref{cndef})
\begin{align}                                  \notag 
c_2 &= 2d (2d-\lambda).
\end{align}
So we have
\begin{align}                                    \notag 
B_2 (x) &=  \frac{\Pi_2(x)}{\lambda 2d (2d-\lambda)}
= -\frac{\delta_{0x}}{2d-\lambda}.
\end{align}
Since $\psi_2(0)= 2^{-d/2} \varphi_\nu (0)=(4\pi \nu)^{-d/2}$, we obtain
\begin{align}                                  \notag 
\abs{B_2 (x)} &\le \decc\,\psi_2(x),\quad\text{whenever}\quad
\decc \ge \frac{(4\pi \nu)^{d/2}}{2d-\lambda}.
\end{align}

Now we come to the induction step: Fix $m\ge 3$ and assume that
$\abs{B_k(x)}\le
\decc
\psi_k(x)$ for all $2\le k< m$. Define the truncated sequence
$(\bbb_n)_{n\ge2}$ by
\begin{align}                                 \notag 
\abs{\bbb_n (x)} &\deff \begin{cases}
B_n(x), &\text{if } \abs{B_n(x)}\le \decc \psi_n(x),\\
\decc \psi_n (x), &\text{if } \abs{B_n(x)}> \decc \psi_n(x).
\end{cases}
\end{align}

This sequence satisfies the hypothesis of Theorem \ref{thm-localper}.
Thus we obtain a sequence $(\baa_n)_{n\in\No}$ of measures with
\begin{align}                                             \label{A'local}
\Babs{\frac{\baa_n (x)}{\ba_n} - 2 \varphi_{n \bd}(x)}
&\le \const  \bigl[n^{-1/2} \varphi_{n\nu}(x) + n^{-d/2} \sum_{k=1}^{n/2}k
\varphi_{k\nu}(x)\bigr],
\end{align}
whenever $n$ has the same parity as $\norme{x}$.
As long as $\lambda$ is small enough, the positive constants $\const$ and
$\nu$
do not depend on $\decc$.

Defining $\bcc_n\deff \bm^n \baa_n$ and using \eqref{A'local} as well as the
fact that $\bd\le \nu$ and both are of comparable size,
we have
\begin{align}                                     \notag 
\bcc_n (x)
&\le \const \bc_n \varphi_{n\nu}(x)
\le \const (\Bar{\alpha}+ \const n^{-1/2}) \bm^n \varphi_{n\nu}(x)
\le \conste  \bm^n \varphi_{n\nu}(x),
\end{align}
where $\conste$ is a positive constant that we fix for the rest of the proof.
Since $\bbb_n$
equals
$B_n$ for all
$n< m$, we also have
$\bcc_n = C_n$ for $n< m$. This can easily be seen by induction.

Now we consider $\Pi_m$. Recall from Lemma
\ref{lem-Pimbound} that
\begin{equation*}
\abs{\Pi_m(x)}\le \lambda \const \conste  \bm^m
m^{-d/2}\sum_{k=1}^{m/2}k^{1-d/2}\,\varphi_{k\nu}(x),
\end{equation*}
and therefore
\begin{align}                                            \label{applystep}
\abs{B_m(x)}&\le
\const \conste  \frac{\bm^m}{c_m}
m^{-d/2}\sum_{k=1}^{m/2}k^{1-d/2}\,\varphi_{k\nu}(x).
\end{align}
It remains to show that $\bm^m/c_m$ is bounded. We know from the lace
expansion formula that
\begin{align}                                            \label{applystep1}
c_m&=2d \bc_{m-1} + \sum_{k=2}^m \pi_k \bc_{m-k}\notag\\ &= 2d
\bm^{m-1}\ba_{m-1} + \sum_{k=2}^m \pi_k \bm^{m-k}\ba_{m-k},
\end{align}
where $\pi_k$ as usual denotes the total mass of $\Pi_k$.
Now we apply Lemma \ref{lem-Pimbound} to the $\pi_k$. We obtain $\abs{\pi_k}
\le
\lambda \const \conste  \bm^k k^{-d/2}$ and insert this in equation
\eqref{applystep1}. This leads to
\begin{align}                                      \notag 
\frac{c_m}{\bm^m}
&\ge d \bm^{-1} - \lambda \const \conste
\ge \const, \text{ if $\lambda=\lambda(d,\conste )$ small enough.}
\end{align}
Using this in equation \eqref{applystep} yields
\begin{align}                                  \notag 
\abs{B_m(x)}&\le
\const \conste m^{-d/2}\sum_{k=1}^{m/2}k^{1-d/2}\,\varphi_{k\nu}(x).
\end{align}
By choosing $\decc$ large enough we see that the lemma follows.
\end{proof}

\begin{proof}[Proof of Theorem \ref{thm-main}.]
Using Lemma \ref{lem-behave}, the second part of the theorem follows directly
from Theorem \ref{thm-localper}. In view of \eqref{bmdecay}, the
first part is a
consequence of   Corollary \ref{cor:speed}.
\end{proof}

\subsection{Identification of the Involved Parameters}

For the convenience of the reader we give a little survey of the parameters
involved in Theorem \ref{thm-main}. The following formulas were first
derived by Brydges and Spencer \cite{BS} (see also \cite{MS} and \cite{NLE}).

The connective constant $\mu$ satisfies the identity
\begin{align}
2d\mu^{-1} = 1-\sum_{m=2}^\infty \pi_m \,\mu^{-m}.\notag
\intertext{The limit $\alpha$ of the mass constants $a_n$ is given by}
\alpha^{-1} = 1+\sum_{m=2}^\infty (m-1) \pi_m \,\mu^{-m},\notag
\intertext{and for the diffusion constant $\delta$ we have}
\delta = \frac{1-\sum_{m=2}^\infty (\pi_m - \underline{\pi}_m) \,\mu^{-m}}%
{d(1+\sum_{m=2}^\infty (m-1)\pi_m \,\mu^{-m})}.        \label{deltaid}
\end{align}

These formulas are obtained by substituting $b_m = \pi_m/(\lambda c_m)$ and
$a_m = \mu^{-m} c_m$ into \eqref{mu}, \eqref{alpha1} and \eqref{deltadef},
respectively.

\appendix
\section{A LCLT and Discretization Estimates}
\label{chapa-normdens}

In order to make the fixed point argument in section \ref{chap-local}, we
need
good local approximations of a general symmetric random walk and quite
specific discretization estimates of
$d$-dimensional normal densities.
We state and prove these in the
         following lemmas. Lemma \ref{lem-approx} is a local central
limit theorem. It controls the pointwise distance between a symmetric
random walk and the appropriate density. This result is proven by
standard large deviation techniques.

Lemma
\ref{lem-A} yields approximations of the discrete folding of a normal
density with a symmetric signed measure on $\Zd$. Analogously, Lemma
\ref{lem-B} approximates the density itself in the variance variable. Both
results are easily calculated by using Taylor expansions.

In addition we give two estimates concerning the discretization of normal
densities. The first one, stated in Lemma \ref{lem-discsum}, gives a simple
bound for the total mass and the \ane{second moment} of a discretized normal
density. The second one compares the
discrete folding of two normal densities with their continuous folding. This
is the content of Lemma \ref{lem-faltung}.

As in the whole paper, $\varphi_\eta$ denotes the density of the
centered normal distribution on
$\Rd$ with covariance matrix $\eta\cdot \Id_d$, that is,
\begin{equation*}
\varphi_\eta (x) = (2\pi\eta)^{-d/2} \,\exp\bigl(-\frac{x^2}{2\eta}\bigr).
\end{equation*}

>\begin{lem}                                             \label{lem-approx}
Let $\AG $ be the single step distribution of an aperiodic, nondegenerate,
symmetric random walk on $\Zd$ with bounded steplength, that is,
$\AG (x)=0$ for all
$\abs{x}>\step$, where $\step $ is fixed.
Let $\mathcal{E} = \eta\cdot
\Id_d$ denote the covariance matrix of $\AG $ and assume $\eta\ge 1/(2d)$.
Then there exist a polynomial $P_4$ of degree four and positive constants
$\const$ and $\nu'$, such that  for all $x$ in $\Zd$ and for all  natural
$n$,
\begin{align}                                      \label{approxlok}
\babs{\AG ^{*n}(x)-\bigl[1+n^{-1}P_4 (x/\sqrt{n}
)\bigr]\varphi_{n\eta}(x)} &\le \const n^{-3/2}\varphi_{n \nu'}(x).
\end{align}
The coefficients of the polynomial depend (rationally) on the
moments of
$\AG $ up to order four, whereas $\const$ and $\nu'$ can be choosen
independently of the specific law of $\AG $, depending only on $d$ and
$\step$.
\end{lem}

\begin{proof}
The proof of the lemma combines standard large deviation properties with the
approximation of $\AG ^{*n}(x)$ obtained by tilting the measure.

Let $Z(t) \deff \sumx \exp(t\cdot x)\AG (x)$ and  $I(\xi)\deff
\sup_{t\in\Rd}\{t\cdot \xi - \log Z(t)\}$. Standard large
deviation theory (see for example \cite{E}) yields a large deviation principle
with entropy function $I$ for the laws of $\AG ^{*n}(x/n)$.
Let $S_\AG$ denote the convex closure of the set of points with nonzero
$\AG$ measure. Then  $I$ is convex on $\Rd$ and even strictly convex on
$\inter S_\AG $, that is, the interior of $S_\AG$. Outside
$S_\AG $, $I$ equals $+ \infty$.

The function $t\mapsto
\nabla\log Z(t)$ is an analytic diffeomorphism from $\Rd$ onto $\inter
S_\AG $ (for a proof see \cite{E}, page 261).
Therefore, for any
$\xi\in \inter S_\AG $, there exists a unique $t_\xi\in\Rd$ with $\nabla\log
Z(t_\xi)=\xi$. Clearly $\nabla \log Z(0)=0$ and $\nabla^2 \log
Z(0)=\mathcal{E}$. For
$\xi \in \inter S_\AG $, we have $I(\xi)= t_\xi\cdot
\xi - \log Z(t_\xi)$. Evidently, $I(0)=0$.
Because of symmetry, the odd partial derivatives of $I$ vanish at zero.
A simple computation yields $\nabla^2 I(0)= \mathcal{E}^{-1}$, and the
fourth derivatives at zero depend only on the second and fourth moments
of $\AG $.

Now denote by $\AG _t$ for $t \in \Rd$ the tilted measure
\begin{equation*}
\AG _t (x) \deff \frac{\AG (x) \exp(t\cdot x)}{Z(t)}.
\end{equation*}
Using this, we see that for $\xi\deff x/n \in \inter S_\AG $, we can write
\begin{equation}                                            \label{Theta-n}
\AG ^{*n} (x) = \exp(-n I(\xi)) \cdot \AG _{t_\xi}^{*n}(x).
\end{equation}

{\em Case $\abs{\xi}\le n^{-5/12}$:}

Since $\AG $ is symmetric and nondegenerate (remember $\eta\ge 1/(2d)$), the
boundary of $S_\AG $ is bounded away from zero. The covariance matrix
$\mathcal{E}_\xi$ of $\AG _{t_\xi}$ is depending
analytically on $\xi$ with $\mathcal{E}_0=\mathcal{E}= \eta\cdot \Id_d$.
It follows that the set $R_\AG$ of all $\xi$ such that
\begin{itemize}
\item
$\xi\in \inter S_\AG$ and
\item
the smallest eigenvalue of $\mathcal{E}_\xi$ is greater or equal $\eta/2$,
\end{itemize}
is a compact neighbourhood of zero.
Thus for $\abs{\xi}\le n^{-5/12}$, we have $\xi\in R_\AG$ for almost all
$n$. It is sufficient to prove the estimate for these $\xi$, since we can
cover the finite number of remaining cases by choosing $\const$ large
enough.

So let $\xi= x/n \in R_\AG $ with $\abs{\xi}\le n^{-5/12}$.
To estimate the first factor in \eqref{Theta-n}, we use
   Taylor expansion for $I$ at zero. We obtain
\begin{align}                                            \label{exp-nI}
\exp(-n I(\xi)) &= \exp(-\frac{x^2}{2n\eta})\cdot
\bigl[1 - n T^{(4)}(\xi) + n O(\xi^6)\bigr]\notag\\
& =\exp(-\frac{x^2}{2n\eta})\cdot
\bigl[1 - n^{-1} T^{(4)}(x/\sqrt{n}) + O(n^{-3/2})\bigr],
\end{align}
where $T^{(k)}$ denotes a polynomial containing $k$th order terms only.
The coefficients of the polynomial are rational functions of the
moments of $\AG $ up to order four.

As a next step we will estimate the second factor in \eqref{Theta-n}, $\AG
_{t_\xi}^{*n}(x)$, using a local central limit theorem out of
Bhattacharya/Rao \cite{BR}. Corollary 22.3 in \cite{BR} asserts
\begin{equation}                                            \label{Theta-t1}
\babs{\AG _{t_\xi}^{*n}(x) -
n^{-d/2}\sum_{r=0}^3 n^{-r/2}Q_r((x-n\xi)/n)\Bigr)}=
o(n^{-(d+3)/2}),
\end{equation}
where $Q_r$ are the so called Edgeworth polynomials. They are
formal polynomials, consisting of partial derivatives of the normal
density in
$\Rd$ with mean zero and covariance matrix $\mathcal{E}_\xi$ (to keep the
notation simple we suppress the $\xi$-dependance of $Q_r$).

For our aims we need the constant implicit in the right hand side of
\eqref{Theta-t1} to be independent of $\xi$, which is a priori not
guaranteed by
\cite{BR}. Calculating the constants in the proof of \cite{BR}, Corollary
22.3, however, shows that they can be choosen such that they only depend
on the maximal steplength $\step$ of
$G_{t_\xi}$ and on a lower bound for the smallest eigenvalue of
$\mathcal{E}_\xi$ on the other. Therefore the error estimate in
\eqref{Theta-t1} is uniform on the compact set $R_\AG$.

Now we come back to the Edgeworth polynomials. In $Q_r$, only derivatives of
order $r+2$, $r+4$,\dots, $3r$ appear (see \cite{BR}, Lemma 7.1). The
coefficients of
$Q_r$ depend on the moments of $\AG _{t_\xi}$ up to order $r+2$. Since $\xi
= x/n$, there is only $Q_r(0)$ appearing in \eqref{Theta-t1}.

$Q_1$ and $Q_3$ vanish at zero, because the odd derivatives
of centered normal densities do so. $Q_0$ is the centered normal
density with covariance matrix $\mathcal{E}_\xi$ itself, so Taylor
expansion yields
$Q_0(0) = (2\pi \eta)^{-d/2}+ T^{(2)}(\xi) + O(\xi^4)$. In the
Taylor expansion of
$Q_2$, the odd terms vanish likewise, and we obtain
$Q_2(0) = \const (2\pi \eta)^{-d/2}  + O(\xi^2)$, where the
constant $\const$ and the error term depend only on the moments of
$\AG $ up to order four. Therefore \eqref{Theta-t1}
simplifies to
\begin{align}                                         \label{Theta-t2}
\AG _{t_\xi}^{*n}(x) &=
(2\pi n \eta)^{-d/2} [ 1+ T^{(2)}(\xi) + O(\xi^4) + n^{-1}T^{(0)}(\xi) +
O(n^{-1}\xi^2) + o(n^{-3/2})]\notag\\
&=
(2\pi n \eta)^{-d/2} [ 1+ n^{-1}T^{(2)}(x/\sqrt{n}) + n^{-1}T^{(0)}
(x/\sqrt{n})
+ O(n^{-3/2})].
\end{align}

Inserting \eqref{exp-nI} and \eqref{Theta-t2} in \eqref{Theta-n} yields
\begin{align}                                  \notag 
\AG ^{*n}(x) &=
\varphi_{n\eta}(x)\cdot
\bigl[1 + n^{-1} P_4(x/\sqrt{n}) + O(n^{-3/2})\bigr],
\end{align}
where $P_4$ is a polynomial of degree four with even order terms only. This
yields the desired estimate whenever $\nu'\ge \eta$.

{\em Case $\abs{\xi}\ge n^{-5/12}$:}

For \ane{big} $x$ we estimate $\AG ^{*n}(x)$ and
$[1 + n^{-1} P_4(x/\sqrt{n})] \varphi_{n\eta}(x)$ separately.
Since for fixed natural $k$ we have $\abs{x}^k \exp (-x^2)\le \const
\exp(-x^2/\sqrt{2})$ for all $x\in\Zd$, the latter is bounded by $\const
\varphi_{\sqrt{2}n\eta}(x)$. Using
$n^{3/2}\le (x/\sqrt{n})^{18}$ in addition, this yields
\begin{equation}                                              \label{pol}
\abs{1 + n^{-1} P_4(x/\sqrt{n})} \varphi_{n\eta}(x)\le \const n^{-3/2}
\varphi_{n\nu'}(x)
\end{equation}
for $\nu'\ge 2\eta$ and for all $n$.

Now consider $\AG ^{*n}(x)$. If $\xi\notin S_\AG $, $\AG ^{*n}(x)$
equals zero.
If $\xi\in S_\AG $, we bound $I(\xi)$ away from zero
with $\theta \,\xi^2$, using Taylor
expansion: We can do this in a neighborhood of zero with the constant
$\theta= 1/(2\eta)$. Since $S_\AG $ is compact (we even have $1/\sqrt{d}\le
\abs{z}\le \step $  for all $z$ in the boundary of $S_\AG $) and $I$ convex,
we
have a nonzero minimum of $I$ on $\partial S_\AG $ and hence we find a
constant
$\theta>0$ such that $I(z)\ge \theta z^2$ for all $z$.
Since $\AG$ is symmetric in each coordinate and rotationally
invariant, $\theta$
can be choosen depending only on $d$ and the range $\step $, but not on the
specific law of $\AG $.

If $\xi\in\inter S_\AG $, we use \eqref{Theta-n} to obtain
$\AG ^{*n} (x) \le \exp(-\frac{x^2}{n\theta})$. The arguments leading to
\eqref{pol} yield the desired estimate for $\nu'\ge 2\theta$.

If $\xi$ lies in the boundary $\partial S_\AG $, use the large deviation
principle to obtain
\begin{align*}                                           
\AG ^{*n}(x) &\le \AG ^{*n}(n\cdot\partial S_\AG )
&\le \exp(-n \inf_{z\in\partial S_\AG }I(z)+n\varepsilon)
&\le \exp(-n \theta\inf_{z\in\partial S_\AG }z^2+n\varepsilon)
\end{align*}
for $n$ large enough. Using  $\inf_{z\in\partial S_\AG }z^2\ge 1/d \ge
\xi^2/(d\step^2)$ and choosing $\varepsilon$ small enough we obtain
$\AG ^{*n} (x) \le \exp(-\frac{x^2}{2d\step^2n\theta})$. The rest of the
argument proceeds as before.

Combining the different bounds and using $\eta\le \step^2/d$ we see that
\eqref{approxlok} holds whenever
$\nu'\ge
\max \,\{4d\step^2\theta, 2\step^2/d\}$.
\end{proof}

\begin{lem}                                                \label{lem-A}
Let $\eta>1/(2d)$ and $\AG \in\mathcal{M}$ (recall the definition at the
beginning of section \ref{chap-local}). Then we have for all
$x\in
\Zd$:
\begin{subequations}\label{A}
\begin{equation}                                                 \label{A1}
\AG *\varphi_{n\eta} (x)= \ag\,\varphi_{n\eta}(x) +
\frac{\ulag}{2dn\eta}\Bigl[\frac{x^2}{n\eta}-d\Bigr]\varphi_{n\eta}(x) +
R_4^\AG (n,x),
\end{equation}
where $R_4^\AG (n,x) = \frac{1}{6}\int_0^1 \!\!ds \,\,(1-s)^3
\sum_{z\in\Zd} \AG (z) D_x^4 \varphi_{n\eta}(x-sz)(z,z,z,z)$.
If $P_{2j}$ is a fixed polynomial of degree $2j$ for some $j \in\No$, then
\begin{equation}                                                 \label{A2}
\AG* [P_{2j}(x/\sqrt{n})\varphi_{n\eta}] (x) =
\ag\,[P_{2j}(x/\sqrt{n})\varphi_{n\eta}](x) + R_2^\AG (n,x;P_{2j}),
\end{equation}
\end{subequations}
where
$R_2^\AG (n,x;P_{2j}) = \int_0^1 \!\!ds \,\,(1-s) \sum_{z\in\Zd}
\AG (z) D_x^2 [P_{2j}(./\sqrt{n})\varphi_{n\eta}](x-sz)(z,z)$.
As local estimates for the error terms we get the following versions,
depending on
$\AG $: If there is a constant $\step$ such that $\AG (x)=0$ for all
$\abs{x}>\step$, we can estimate the error terms by
\begin{equation}
\begin{split}
\abs{R_{4}^\AG (n,x)} &\le \const(d,\eta,\step)\;n^{-2}
\label{R}
\,\,\sum_{z\in\Zd} z^4\, \abs{\AG (z)}\;
\varphi_{2n\eta}(x)\quad\text{and}\\
\abs{R_{2}^\AG (n,x;P_{2j})} &\le \const(d,\eta,j,\step)\;n^{-1}
\,\,\sum_{z\in\Zd} z^2\, \abs{\AG (z)}\;\varphi_{2n\eta}(x).
\end{split}
\end{equation}
If no such $\step$ exists, we still have the
estimates
\begin{equation}
\begin{split}
\abs{R_4^\AG (n,x)} &\le \const(d,\eta)\; n^{-2}
\label{R'}
\,\sum_{z\in\Zd} z^4\, \abs{\AG (z)}\,
\int_0^1 \!\!ds\,\,\varphi_{\sqrt{2}n\eta}(x-sz)
\quad\text{and}\\
\abs{R_2^\AG (n,x;P_{2j})} &\le \const(d,\eta,j)\; n^{-1}
\,\sum_{z\in\Zd} z^2\, \abs{\AG (z)}\,
\int_0^1 \!\!ds\,\,\varphi_{\sqrt{2}n\eta}(x-sz).
\end{split}
\end{equation}
\end{lem}

\begin{proof}
Using Taylor expansion and the symmetry of $\AG $ we obtain
\begin{align*}
\AG *\varphi_{n\eta}(x) &= \sum_{z\in\Zd}\AG (z)\varphi_{n\eta}(x-z) \notag\\
&=\ag \varphi_{n\eta}(x) + \frac{\ulag }{2d}\Delta_x\varphi_{n\eta}(x)
+ R_4^\AG (n,x),
\end{align*}
which is \eqref{A1} after inserting $\Delta_x\varphi_{n\eta}(x) =
\bigl[\frac{x^2}{n^2\eta^2}-\frac{d}{n\eta}\bigr]\varphi_{n\eta}(x)$.
Analogously, first order Taylor approximation leads to
equation \eqref{A2}.

Now we come to the proof of the error estimates.
We write $P_{(k)}$ to denote some polynomial of order $k$.
The forth partial
derivatives of $\varphi_{n\eta}$ are functions of the form
$n^{-2}P_{(4)}(./\sqrt{n})\varphi_{n\eta}$ and therefore bounded by
$\const(d,\eta) n^{-2}\varphi_{\sqrt{2}n\eta}$. Similarly, the second
partial derivatives of
$P_{2j}(./\sqrt{n})\varphi_{n\eta}$ are of the form
$n^{-1}P_{(2j+2)}(./\sqrt{n})\varphi_{n\eta}$ and bounded by
$\const(d,\eta,j)\; n^{-2}\varphi_{\sqrt{2}n\eta}$. This implies \eqref{R'}.

To prove \eqref{R}, we use the fact that for
$z$ with $\abs{z}\le \step$ we have $\varphi_{\sqrt{2}n\eta}(x-z)\le \const
\varphi_{2n\eta}(x)$
with $\const$ depending on $d$, $\eta$ and $\step$, but not on $n$.
\end{proof}

\begin{lem}                                                   \label{lem-B}
Let $\eta>0$. Then for $x\in\Zd$, $n\in\N$ and $k\in(0,n)$:
\begin{subequations}\label{B}
\begin{equation}                                                 \label{B1}
\varphi_{(n-k)\eta}(x) = \varphi_{n\eta}(x) -
\frac{k}{2n}\Bigl[\frac{x^2}{n\eta}-d\Bigr]\varphi_{n\eta}(x) +
S_2^k(n,x),
\end{equation}
where $S_2^k(n,x) = \frac{1}{2}k^2 \frac{\partial^2}{\partial
\theta^2}\varphi_{\theta \eta}(x)$ for a $\theta\in(n-k,n)$.
On the other hand, for a fixed polynomial $P_{2j}$ of degree $2j$,
$j\in\No$:
\begin{equation}                                                 \label{B2}
[(n-k)^{-1}P_{2j}(x/\sqrt{n-k})\varphi_{(n-k)\eta}](x) =
n^{-1} P_{2j}(x/\sqrt{n})\varphi_{n\eta}(x) + S_1^k(n,x;n^{-1}P_{2j}),
\end{equation}
\end{subequations}
where $S_1^k(n,x;n^{-1}P_{2j}) = - k \frac{\partial}{\partial
\theta}P_{2j}(x/\sqrt{\theta})\varphi_{\theta \eta}(x)$ for some
$\theta\in(n-k,n)$.
Estimates for the error terms are given by
\begin{equation}
\begin{split}
\abs{S_2^k(n,x)} &\le \const(d,\eta)\;\Bigl(\frac{k}{n-k}\Bigr)^{2}
\label{S}
\Bigl(\frac{n}{n-k}\Bigr)^{d/2}\varphi_{\sqrt{2}n\eta}(x)\quad\text{and}\\
\abs{S_1^k(n,x;n^{-1}P_{2j})} &\le \const(d,\eta,j)\;\frac{k}{(n-k)^2}
\Bigl(\frac{n}{n-k}\Bigr)^{d/2}\varphi_{\sqrt{2}n\eta}(x).
\end{split}
\end{equation}

\end{lem}

\begin{proof}
Here we use one dimensional Taylor expansion for $\varphi_{n\eta}(x)$ as a
function in $n$ to write
\begin{align*}
\varphi_{(n-k)\eta}(x) &= \varphi_{n\eta}(x) -
k \frac{\partial}{\partial n}\varphi_{n\eta}(x) + S_2^k(n,x),
\end{align*}
which implies \eqref{B1} by using $\frac{\partial}{\partial
n}\varphi_{n\eta}(x) =
\bigl[\frac{x^2}{2\eta n^2}-\frac{d}{2n}\bigr]\varphi_{n\eta}(x)$. Keeping
only
the constant term of the Taylor approximation leads to
\eqref{B2}.

To prove the error estimates, we first observe that the second derivative
(with respect to
$\theta$) of $\varphi_{\theta \eta}(x)$ is of the form
$\theta^{-2}P_{(4)}(x/\sqrt{\theta})\varphi_{\theta \eta}(x)$ and therefore
bounded by
$\const(d,\eta)
\theta^{-2} \varphi_{\sqrt{2}\theta \eta}(x)$, while the first derivative of
$\theta^{-1}P_{2j}(x/\sqrt{\theta})\varphi_{\theta \eta}(x)$ has the form $
\theta^{-2}P_{(2j+2)}(x/\sqrt{\theta})\varphi_{\theta \eta}(x)$ and bounded
by
$\const(d,\eta,j)
\theta^{-2} \varphi_{\sqrt{2}\theta \eta}$. Again, $P_{(k)}$ stands for a
polynomial of order $k$.

Splitting the function $\varphi$
and replacing $\theta$ separately by $n-k$ in the  first factor
and by $n$ in the exponential term leads to \eqref{S}.
\end{proof}

\begin{lem}                                             \label{lem-discsum}
Assume $\eta \ge 1/(2d)$. Then there exists a constant $\const$ depending only
on the dimension $d$, such that
\begin{align}                               \notag 
\sumx \varphi_\eta (x)&\le \const&&
\text{and}
&\sumx x^2 \varphi_\eta (x)&\le \const \eta.
\end{align}
\end{lem}

\begin{proof}
For $d=1$ we have
\begin{align*}
\sum_{x\in\Z} \varphi_\eta (x)
= 2 \sumn \varphi_\eta (n) +  \varphi_\eta(0)
\le 2\int_0^\infty \varphi_\eta (t) dt + \frac{1}{\sqrt{2\pi\eta}}
\le 1 + \frac{1}{\sqrt{2\pi\eta}}\le \const.
\end{align*}
The first inequality for general dimension follows immediately from this
estimate, because
the sum over $\Zd$ of the density values equals the $d$th power of the sum
over
$\Z$ of the values of the one dimensional normal density with  variance
$\eta$.

The second inequality follows from the first one by using
\begin{align*}
&\sum_{x\in\Zd} x^2 \varphi_\eta (x)
\le  \const\eta\,\sum_{x\in\Zd} \varphi_\eta (x),
\end{align*}
which comes from the fact that $(x^2/\eta)\exp(-x^2/(2\eta))$ can
be  bounded  uniformly in $x\in\Zd$ by
$\const\,\exp(-x^2/(4\eta))$.
\end{proof}

\begin{lem}                                             \label{lem-faltung}
Let $\eta,\theta\geq 1/(2d).$ Then for $x\in \Zd$%
\begin{equation*}
\sum_{y\in \Zd}\varphi _{\eta}(y)\varphi _{\theta}(x-y)\leq \const\varphi
_{\eta+\theta}(x).
\end{equation*}
\end{lem}

\begin{proof}
Let $I^d\deff\bigl[-\frac{1}{2},\frac{1}{2}\bigr]^{d}$,
and denote by $y+I^d$ the shifted cube. Then
\begin{align*}
\varphi _{\eta+\theta}(x) &=\sum_{y\in \Zd}\int_{y+I^d}\varphi
_{\eta}(t)\varphi _{\theta}(x-t)\,dt \\
&=\frac{1}{( 2\pi \eta) ^{d/2} ( 2\pi \theta) ^{d/2}}%
\sum_{y\in \Zd}\int_{I^d}\exp \bigl[ -\frac{1}{2\eta}(y+t)^{2}
-\frac{1}{2\theta}(x-y-t)^{2}\bigr] \,dt \\ &\ge \frac{1}{\bigl(
2\pi \eta\bigr) ^{d/2}\bigl( 2\pi \theta\bigr) ^{d/2}}%
\sum_{y\in \Zd}\exp \bigl[ -\int_{I^d}\bigl(
\frac{1}{2\eta} (y+t)^{2}+\frac{1}{2\theta}(x-y-t)^{2}\bigr)
\,dt\bigr]
\end{align*}
by Jensen's inequality. Note that for $t\in I^d$ we have $t^2\le d/4$,
so that
\begin{align*}
\int_{I^d}\Bigl( \frac{1}{2\eta}(y+t)^{2}+\frac{1}{2\theta}(x-y-t)^{2}
\Bigr) \,dt &\le \frac{1}{2\eta}y^{2}+\frac{1}{2\theta}
(x-y)^{2}+\frac{d}{8\eta}+\frac{d}{8\theta} \\
&\leq\frac{1}{2\eta}y^{2}+\frac{1}{2\theta}(x-y)^{2}+\const,
\end{align*}
because we assumed that $\eta,\theta\geq 1/(2d)$. Therefore
\begin{flalign*}
{}&&\varphi_{\eta+\theta}(x)&\ge \mathrm{e}^{-\const}\sum_{y\in
\Zd}\varphi _{\eta}(y)\varphi _{\theta}(x-y).&&\qed
\end{flalign*}
\renewcommand{\qedsymbol}{}\end{proof}
\section{The Lace Expansion}                              \label{chap-le}

This section contains standard material on the lace expansion in the first
subsection and bounds for the lace expansion terms in the second one. The
lace
expansion was introduced by Brydges and Spencer in \cite{BS} and discussed in
detail by Madras and Slade in \cite{MS}. The following overview consists of
the
minimum necessary to make this thesis selfcontained. The first part is taken
more or less literally from van der Hofstad, den Hollander and Slade
\cite{NLE}.

\subsection{Definition of the Lace Functions}

In this section we define the Lace Functions $\Pi_m$ and prove the recursion
formula \eqref{recursion}, that is
\begin{equation*}
       C_n = 2d(D*C_{n-1}) + \sum_{m=2}^n \Pi_m*C_{n-m}.
\end{equation*}
This requires the introduction of the following
standard terminology.
Given an interval $I=[a,b]\subset \Z$ of integers with $0\le a\le b$, we refer
to a pair
$\{s,t\}$ ($s<t$) of elements of $I$ as an {\em edge}. To abbreviate
the notation, we write $st$ for $\{s,t\}$. A set of edges is called a {\em
graph}. A graph $\Gamma$ on $[a,b]$ is said to be {\em connected} if both $a$
and $b$ are endpoints of edges in $\Gamma$ and if, in addition, for any
$c\in[a,b]$ there is an edge $st\in\Gamma$ such that $s<c<t$. The set of all
graphs on
$[a,b]$ is denoted $\mathcal{B}[a,b]$, and the subset consisting of all
connected graphs is denoted $\mathcal{G}[a,b]$. A {\em lace} is a
minimally connected graph, that is, a connected graph for which the
removal of any edge would result in a disconnected graph. The set of
laces on $[a,b]$ is denoted $\mathcal{L}[a,b]$, and the set of laces on
$[a,b]$ consisting of exactly $N$ edges is denoted
$\mathcal{L}^{(N)}[a,b]$.

Given a connected graph $\Gamma$, the following prescription associates to
$\Gamma$ an unique lace $L_\Gamma$: The lace consists of edges $s_1t_1,
s_2t_2,\dots$, with $t_1$, $s_1$, $t_2$, $s_2$, $\dots$ determined (in that
order) by
\begin{align*}
t_1 &= \max\{t:\, at\in\Gamma\},& s_1 &= a,\\
t_{i+1}&= \max\{t:\, \exists s< t_i\,\text{ such that }st\in\Gamma\},&
s_{i+1}&= \min\{s:\, st_{i+1}\in\Gamma\}.
\end{align*}

Given a lace $L$, the set of all edges $st\notin L$ such that $L_{L \cup
\{st\}}= L$ is denoted $\mathcal{C}(L)$. Edges in $\mathcal{C}(L)$ are
said to be {\em compatible} with $L$.

Recall the definition of $C_n$,
\begin{equation*}
C_n(x) = \sumsx \prod_{0\le s<t\le n} (1-\lambda U_{st}(\omega)),
\end{equation*}
where $U_{st}(\omega) =\delta_{\omega(s),\omega(t)}$. Now we define
for integers
$0\le a<b$
\begin{align}                                        \label{Kabdef}
K[a,b](\omega) \deff \prod_{a\le s<t \le b}(1-\lambda U_{st}(\omega)).
\end{align}
Then we can write
\begin{align}                                         \label{Cnneu}
C_n(x) = \sumsx K[0,n](\omega),
\end{align}
where the sum is over all $n$ step simple random walk paths from $0$ to $x$.
Expanding the product in the definition of $K[a,b](\omega)$, we get%
\begin{align}                                                     \label{Kab}
K[a,b](\omega) =\sum_{\Gamma\in \mathcal{B}[a,b]} \prod_{st\in\Gamma}(-\lambda
U_{st}(\omega)).
\end{align}
We also define an analogous quantity, in which the sum over graphs
is restricted to connected graphs, namely,%
\begin{align}                                         \label{Jabdef}
J[a,b](\omega) \deff\sum_{\Gamma\in \mathcal{G}[a,b]}
\prod_{st\in\Gamma}(-\lambda
U_{st}(\omega)).
\end{align}
This allows us to define the \ane{lace functions}, which are the key
quantities
in the lace expansion,%
\begin{align}                                         \label{Pimdef}
\Pi_m(x) =\sumsx[=m] J[0,m](\omega).
\end{align}
The identity \eqref{recursion} now follows from the following lemma.

\begin{lem}                                             \label{lem-recursion}
For $n\ge 1$,
\begin{align}                                         \label{recurs1}
C_n(x) = \sum_{y:\norme{y}=1} C_{n-1}(x-y) +
\sum_{m=2}^n \sum_{z\in\Zd} \Pi_m(z)\, C_{n-m}(x-z).
\end{align}
\end{lem}

\begin{proof}
It suffices to show that for each path $\omega$ we have (suppressing $\omega$
in the formulas):
\begin{align}                                         \label{recurs2}
K[0,n] = K[1,n] + \sum_{m=2}^n J[0,m] \, K[m,n].
\end{align}
Then \eqref{recurs1} is obtained after insertion of \eqref{recurs2} into
\eqref{Cnneu} followed by factorization of the sum over $\omega$. To prove
\eqref{recurs2}, we note from \eqref{Kab} that the contribution to $K[0,n]$
from all graphs $\Gamma$ for which $0$ is not in an edge is exactly $K[1,n]$.
To resum the contribution from the remaining graphs, we proceed as follows.
When $\Gamma$ does contain an edge ending at $0$, we let $m[\Gamma]$ denote
the largest value of $m$ such that the set of edges in $\Gamma$ with at least
one end in the interval $[0,m]$ forms a connected graph on $[0,m]$. We lose
nothing by taking $m\ge 2$, since $U_{a,a+1}=0$ for all $a$. Then resummation
over graphs on $[m,n]$ gives
\begin{align}                                \notag 
K[0,n] = K[1,n] + \sum_{m=2}^n \sum_{\Gamma\in \mathcal{G}[0,m]} \prod_{st\in
\Gamma}(-\lambda U_{st})\, K[m,n].
\end{align}
With \eqref{Jabdef} this proves \eqref{recurs2}.
\end{proof}

We next rewrite \eqref{Pimdef} in a form that can be used to obtain good
bounds on $\Pi_m(x)$. For this, we begin by partially resumming the right-hand
side of \eqref{Jabdef}, to obtain
\begin{align}                                         \label{Jab}
J[a,b] &= \sum_{L\in\mathcal{L}[a,b]} \sum_{\Gamma:\, L_\Gamma = L}
\prod_{st\in L}(-\lambda U_{st})\, \prod_{s't'\in \Gamma\setminus L}(-\lambda
U_{s't'})\notag\\
&= \sum_{L\in\mathcal{L}[a,b]} \prod_{st\in L}(-\lambda U_{st})\,
\prod_{s't'\in
\mathcal{C}(L)}(1-\lambda U_{s't'}).
\end{align}

For $0\le a<b$, we define $J^{(N)}[a,b]$ to be, up to the factor
$(-\lambda)^{N}$, the contribution to
\eqref{Jab} coming from laces consisting of exactly $N$ edges,
\begin{align}                                \notag 
J^{(N)}[a,b] &\deff
\sum_{L\in\mathcal{L}^{(N)}[a,b]} \prod_{st\in L} U_{st}\,
\prod_{s't'\in \mathcal{C}(L)}(1-\lambda U_{s't'}), \qquad N\ge 1.
\end{align}
Then
\begin{align}                              \notag 
J[a,b] &= \sum_{N=1}^\infty (-\lambda)^N J^{(N)}[a,b],
\end{align}
and by \eqref{Pimdef}
\begin{align}                                         \label{Pim}
\Pi_m(x) &= \sum_{N=1}^\infty (-\lambda)^N \Pi_m^{(N)}(x),
\end{align}
where we define
\begin{align}                                         \label{PimNdef}
\Pi_m^{(N)}(x) &\deff \sumsx[=m] J^{(N)}[0,m](\omega)\notag\\
&= \sumsx[=m] \sum_{L\in\mathcal{L}^{(N)}[0,m]} \prod_{st\in L} U_{st}\,
\prod_{s't'\in \mathcal{C}(L)}(1-\lambda U_{s't'}).
\end{align}

\subsection{Bounds on the Lace Functions}

In this section, we obtain bounds on $\Pi_m(x)$. A first part will provide the
standard bounds in terms of $C_n(x)$, the two-point functions of the
appropriate weakly self avoiding walk after $n$ steps. In the second part we
will obtain specific bounds by assuming a Gaussian decay of the $C_n(x)$. In
this case good bounds result easily from the Cauchy-Schwarz inequality.

There is only one lace on $[0,m]$ consisting of exactly one edge, namely
$\{0m\}$. Therefore we have for $N=1$:
\begin{align}                                         \label{Pim1}
\Pi_m^{(1)}(x) &= \delta_{0x}\sums[=m]{0}{0} \;\;
\prod_{\substack{0\le s'<t'\le m\\ s't'\neq 0m}}(1-\lambda U_{s't'}(\omega))
\notag\\&
\le \delta_{0x}\sum_{y:\,\norme{y}=1}\sums[=m-1]{y}{0} K[1,m](\omega)
&&\text{using }1-U_{0t'}\le 1 \quad\forall \,t' \text{ and }\eqref{Kabdef}
\notag\\&
= \delta_{0x} \;2d D* C_{m-1}(0)&&\text{by }\eqref{Cnneu}.
\end{align}

Now for $N\ge 2$: A walk giving a nonzero contribution to \eqref{PimNdef}
must
intersect itself
$N$ times, to ensure that $U_{st}\neq 0$ for each $st \in L$. For example,
when
$N=7$, the walk must undergo a trajectory of the form
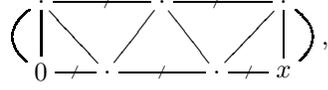
\begin{equation}                                \notag 
\raisebox{0.5 \depth}{%
\xymatrix{%
&.\ar@{-}[d]\ar@{-}[rrrr]&&\strich&&.\ar@{-}[dll]%
\ar@{-}[rrrr]&&\strich&&.\ar@{-}[dll]%
\ar@{-}@/^\krue em/[d]&%
\\
&0\ar@{-}@/^\krue em/[u]\ar@{-}[rr] &\strich& .\ar@{-}[ull]%
\ar@{-}[rrrr] &&\strich&&.\ar@{-}[ull]%
\ar@{-}[rr] &\strich& x\ar@{-}[u]}},
\end{equation}
where the slashed lines denote subwalks that may have length zero, whereas
nonslashed lines denote subwalks containing at least one step.

Using $1-\lambda U_{s't'}\le 1$ in
\eqref{PimNdef} whenever $s'$ and
$t'$ belong to different subwalks, we get an upper bound in which distinct
subwalks no longer interact. However, each subwalk remains weakly
self-avoiding. We thus have
\begin{align}                                            \label{PimNest}
\Pi_m^{(N)}(x)&\le
\sum_{y_i,m_i}C_{m_1}(y_1)C_{m_2}(y_1)C_{m_3}(y_2)C_{m_4}(y_1-y_2)
C_{m_5}(y_3-y_1)\cdots\\
&\hspace{2cm}\cdots
C_{m_{2N-3}}(x-y_4) C_{m_{2N-2}}(y_5-x) C_{m_{2N-1}}(x-y_5),\notag
\end{align}
where the sum is over $y_1,\dots,y_{N-2}\in\Zd$ and over $m_1,\dots,m_{2N-1}$
with $\sum m_i=m$. The $m_i$ are nonnegative integers, and only $m_{3}$,
$m_{5}$, $\dots$, $m_{2N-3}$ can equal zero.

\begin{lem}                                               \label{lem-Pimbound}
Fix $m\ge 2$ and $d\ge 5$.
Assume that for all $x\in \Zd$ and natural $n< m$, we have
\begin{equation}                                  \label{Cnest}
\abs{C_n(x)}\le \conste  \mu^n \varphi_{n\nu}(x)
\end{equation}
with constants $\mu>0$, $\conste \ge 1$ and $\nu\ge 1/(2d)$. Then, for
$\lambda=\lambda(d,\conste )$ small enough, we have
\begin{equation}
\abs{\Pi_m(x)}\le \lambda \const \conste  \mu^m
m^{-d/2}\sum_{k=1}^{m/2}k^{1-d/2}\,\varphi_{k\nu}(x).
\end{equation}
\end{lem}

\begin{proof}
Again we abbreviate the notation by writing
\begin{equation*}
\psi_m (x) \deff m^{-d/2}\sum_{k=1}^{m/2}k^{1-d/2}\varphi_{k \nu}(x).
\end{equation*}
We consider the terms $\Pi_m^{(N)}$ separately.
For $N=1$, equation \eqref{Pim1} yields
\begin{align}                                         \label{Pim1'}
\Pi_m^{(1)}(x) &
\le \delta_{0x} \;2d \conste  \mu^{m}\, D* \varphi_{(m-1)\nu}(0)
\le \const \conste  \mu^{m}\, \psi_m(x),
\end{align}
since $D* \varphi_{(m-1)\nu}(0)\le \const \varphi_{2(m-1)\nu}(0)\le \const
m^{-d/2}\varphi_{\nu}(0)$.

To keep the notation simple, we set $\varphi_0(x)\deff \delta_{0,x}$.
For $N\ge 2$ we define
\begin{equation}                                   \notag 
\begin{split}
P_m^{(N)}(x)&\deff
\sum_{y_i,m_i}\varphi_{m_1\nu}(y_1)\varphi_{m_2\nu}(y_1)\varphi_{m_3\nu}(y_2)
\varphi_{m_4\nu}(y_1-y_2)\varphi_{m_5\nu}(y_3-y_1)\!\cdots\\
&\hspace{2.1cm}\cdots
\varphi_{m_{2N-3}\nu}(x-y_4) \varphi_{m_{2N-2}\nu}(y_5-x)
\varphi_{m_{2N-1}\nu}(x-y_5),
\end{split}
\end{equation}
with the same summation as in \eqref{PimNest} up to the fact that we
allow an additional term for $m_2=0$ (this corresponds to slashing the line
from zero vertically up and is necessary to give the induction below). By
\eqref{Cnest}, we have  for all $n$:
\begin{equation}                             	\label{PimNest1}
\Pi_m^{(N)}\le \conste^{2N-1}\,\mu^{m}\, P_m^{(N)}.
\end{equation}

Now we show by induction that there is a constant $\constz$ depending only on
the dimension $d$ such that
\begin{equation}                                          \label{induction1}
\abs{P_m^{(N)}}\le \constz^N \psi_m.
\end{equation}

\vspace{1ex}
For $N=2$ the lace diagram is three-legged: \quad
\raisebox{0.5 \depth}{%
\xymatrix{%
0\ar@{-}@/_\krue em/[rrrr]\ar@{-}@/^\krue
em/[rrrr]&&\strich&&x\ar@{-}[llll]}}.
\quad We have

\begin{align*}
P_m^{(2)}(x) &=\sum_{\substack{k+l+j= m\\l,j\ge1}}
\varphi_{k\nu}(x)\varphi_{l\nu}(x)\varphi_{j\nu}(x)
= \delta_{0,x}I+J,
\intertext{where}
I&=\sum_{l=1}^{m-1}
\varphi_{l\nu}(0)\varphi_{(m-l)\nu}(0)
\le \const (2\pi\nu)^{-d}\sum_{l=1}^{m-1} l^{-d/2}(m-l)^{-d/2}\notag\\
&\le \const \underbrace{(2\pi\nu)^{-d/2}}_{=\varphi_\nu(0)}m^{-d/2}\le \const
\psi_m(0)\\
\text{and}\quad
J&\le \const \sum_{\substack{k+l+j= m\\1\le k\le l\le j}}
\varphi_{k\nu}(x)\, l^{-d/2}\, m^{-d/2}\notag\\
&\le \const m^{-d/2}\sum_{k=1}^{m/3} \varphi_{k\nu}(x) \sum_{l=k}^{m}
l^{-d/2}
\le \const \psi_m(x).
\end{align*}
So together we obtain
\begin{align}                                          \label{anker}
P_m^{(2)}(x)\le \const\psi_m(x).
\end{align}

Now we come to the induction step. For $N\ge 3$ we will reduce
$P_m^{(N)}$ to
$P_m^{(N-1)}$ by merging four subwalks in the lace into two. The following
figure illustrates this process:
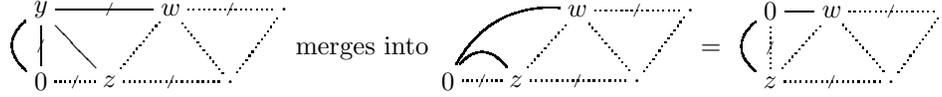
\xymatrixrowsep{.2em}
\xymatrixcolsep{.25em}
\begin{align*}
\raisebox{0.5 \depth}{%
\xymatrix{%
&y\ar@{-}@/_\krue em/[dd]\ar@{-}[dd]&&\strich&&w\ar@{.}[ddll]\ar@{-}[llll]%
&&\strich&&.\ar@{.}[ddll]\ar@{.}[llll]\\%
&\strich\\%
&0 &\strich&z\ar@{.}[ll]\ar@{-}[uull]&&\strich&&.\ar@{.}[llll]\ar@{.}[uull]}}
%
\text{merges into}\!\!\!\!
\xymatrixrowsep{1.5em}
%
\raisebox{0.5 \depth}{%
\xymatrix{%
&&&&&w\ar@{-}@/_\krue em/[dllll]\ar@{.}[dll]%
&&\strich&&.\ar@{.}[dll]\ar@{.}[llll]\\%
&0 \ar@{-}@/^\krue
em/[rr]&\strich&z\ar@{.}[ll]&&\strich&&.\ar@{.}[llll]\ar@{.}[ull]}}
%
\!\!=\,
\xymatrixrowsep{.2em}
%
\raisebox{0.5 \depth}{%
\xymatrix{%
&0\ar@{-}@/_\krue em/[dd]\ar@{.}[dd]&&w\ar@{.}[ddll]\ar@{-}[ll]%
&&\strich&&.\ar@{.}[ddll]\ar@{.}[llll]\\%
&\strich\\%
&z&&\strich&&.\ar@{.}[llll]\ar@{.}[uull]}}
\end{align*}

\xymatrixrowsep{1.5em}
\xymatrixcolsep{0.25em}

We use Cauchy-Schwarz to
obtain
\begin{align}                                            \label{cauchyschwarz}
&\sum_{\substack{u_1+u_2=u\\u_1\ge 1}}\;
\sum_{\substack{t_1+t_2=t\\t_2\ge 1}}\;\sum_{y\in\Zd}
\varphi_{u_1 \nu}(y) \varphi_{u_2\nu}(w-y)
\varphi_{t_1 \nu}(y) \varphi_{t_2\nu}(y-z)\notag\\
&\quad\le \sum_{\substack{u_1+u_2=u\\t_1+t_2=t\\u_1,t_2\ge 1}}
\Bigl[\sum_{y\in\Zd}\varphi^2_{u_1 \nu}(y) \varphi^2_{u_2\nu}(w-y)\Bigr]^{1/2}
\,\Bigl[\sum_{y\in\Zd}\varphi^2_{t_1 \nu}(y)
\varphi^2_{t_2\nu}(y-z)\Bigr]^{1/2}.
\end{align}
Note that for $u_2\ge 1$ we have
\begin{align*}                                            
\bigl[\varphi^2_{u_1 \nu}* \varphi^2_{u_2\nu}\bigr]^{1/2}
&\le \const \nu^{-d/2}(u_1u_2)^{-d/4}
[\varphi_{u\nu/2}]^{1/2}\le \const (u_1u_2)^{-d/4}u^{d/4}
\varphi_{u\nu},
\end{align*}
and therefore (recall $d\ge 5$)
\begin{align*}                                            
\sum_{\substack{u_1+u_2=u\\u_1\ge 1}} \bigl[\varphi^2_{u_1 \nu}*
\varphi^2_{u_2\nu}\bigr]^{1/2}
       &\le \bigl(1+u^{d/4}
\sum_{u_1=1}^{u-1}u_1^{-d/4}(u-u_1)^{-d/4}\bigr)\,\varphi_{u\nu}\le \const
\varphi_{u\nu}.
\end{align*}
Inserting this into \eqref{cauchyschwarz} yields
\begin{align}                                            \label{step1}
\sum&
\varphi_{u_1 \nu}(y) \varphi_{u_2\nu}(w-y)
\varphi_{t_1 \nu}(y) \varphi_{t_2\nu}(y-z)
\le  \const \varphi_{u\nu}(w)\; \varphi_{t\nu}(z).
\end{align}

Using \eqref{step1} with $y_1$, $y_2$, $y_3$, $m_1$, $m_5$, $m_2$ and $m_4$ in
place of $y$, $z$, $w$,
$u_1$, $u_2$, $t_1$ and $t_2$, respectively, gives
\begin{align}                                            \label{step2}
P_m^{(N)}(x)&=\sum_{y_i,m_i} \varphi_{m_2\nu}(y_1) \varphi_{m_4\nu}(y_1-y_2)
\varphi_{m_1\nu}(y_1) \varphi_{m_5\nu}(y_3-y_1)
\varphi_{m_3\nu}(y_2)\cdots
\notag\\
&\le \const
\sum \varphi_{(m_2+m_4)\nu}(y_2) \varphi_{m_3\nu}(y_2)
\varphi_{(m_1+m_5)\nu}(y_3)\cdots\notag\\
&\le \const P_m^{(N-1)}(x).
\end{align}
Now we choose $\constz$ to be the maximum of the constants appearing in
\eqref{anker}  and \eqref{step2}. We obtain \eqref{induction1}. Now
\eqref{Pim1'} and \eqref{PimNest1} together imply
\begin{align}                                            \label{PimNest2}
\Pi_m^{(N)}\le \const \constz^N \, \conste^{2N-1} \, \mu^{m}\,\psi_m.
\end{align}
The lemma now follows by inserting \eqref{PimNest2}
into \eqref{Pim} and choosing $\lambda$ small enough.
\end{proof}


\end{document}